\documentclass[a4paper,12pt,english]{article}
\usepackage[paper=a4paper,left=30mm,right=20mm,top=25mm,bottom=30mm]{geometry}
\usepackage{amsmath}
\usepackage{amsfonts}
\usepackage{amssymb}
\usepackage[all]{xypic}
\usepackage{graphics}
\usepackage{mathtools}
\usepackage{babel}
\usepackage{here}
\usepackage{verbatim}

\def\a{\alpha}
\def\b{\beta}

\def\e{\varepsilon}

\def\vpi{\varpi}
\def\g{\gamma}
\def\h{\eta}

\def\k{\kappa}
\def\l{\lambda}
\def\o{\omega}

\def\s{\sigma}

\def\D{\Delta}

\def\G{\Gamma}

\def\bR{{\mathbb R}}
\def\bZ{{\mathbb Z}}
\def\bC{{\mathbb C}}

\def\b1{{\rm id}}

\newfont{\goth}{eufm10 scaled \magstep1}

\def\gc{\mbox{\goth c}}
\def\gd{\mbox{\goth d}}
\def\gsl{\mbox{\goth sl}}

\def\gso{\mbox{\goth so}}
\def\ggl{\mbox{\goth gl}}

\def\gg{\mbox{\goth g}}
\def\gh{{\mbox{\goth h}}}
\def\gk{\mbox{\goth k}}
\def\gl{\mbox{\goth l}}
\def\gm{\mbox{\goth m}}
\def\gn{\mbox{\goth n}}

\def\gp{\mbox{\goth p}}

\newfont{\mcal}{eusm10 scaled \magstep1}
\def\ca{\mbox{\mcal A}}

\def\cd{\mbox{\mcal D}}

\def\cf{\mbox{\mcal F}}

\def\cl{\mbox{\mcal L}}

\def\cx{\mbox{\mcal X}}

\newtheorem{Th}{Theorem}
\newtheorem{Rem}[Th]{Remark}
\newtheorem{Prop}[Th]{Proposition}
\newtheorem{Cor}[Th]{Corollary}
\newtheorem{Lem}[Th]{Lemma}
\newtheorem{Def}[Th]{Definition}
\newtheorem{Ex}[Th]{Example }
\def\bex{\begin{Ex}}
\def\eex{\end{Ex}}

\def\bt{\begin{Th}}
\def\et{\end{Th}}
\def\bp{\begin{Prop}}
\def\ep{\end{Prop}}
\def\bc{\begin{Cor}}
\def\ec{\end{Cor}}
\def\bl{\begin{Lem}}
\def\el{\end{Lem}}
\def\bd{\begin{Def}}
\def\ed{\end{Def}}
\def\bex{\begin{Ex}}
\def\eex{\end{Ex}}
\def\br{\begin{Rem}}
\def\er{\end{Rem}}
\def\pf{\noindent{\it Proof:\ }}
\def\qed{\hfill$\square$}
\def\n{\nabla}
\def\op{\oplus}
\def\ot{\otimes}

\def\be{\begin{equation}}
\def\ee{\end{equation}}
\def\ben{\begin{enumerate}}
 \def\een{\end{enumerate}}
\def\ba{\begin{array}{rlll}}
\def\ea{\end{array}}
\def\bea{\begin{eqnarray}}

\def\SR{\rm  sub-Riemannian \,}
\def\eea{\end{eqnarray}}
\def\bean{\begin{eqnarray*}}
\def\eean{\end{eqnarray*}}

\def\ad{{\rm ad}\,}
\def\Ad{{\rm Ad}\,}

\def\ker{\mathrm{ker\;}}

\def\id{\rm id \;}

\def\C1{cohomogeneity one  }

\def\p{\partial}

\begin{document}

\title{Shortest  and Straightest geodesics in  Sub-Riemannian Geometry }
	\author{Dmitri Alekseevsky
\footnote{The  work  was partially  supported   by   the grant     no. 18-00496S of
the Czech Science Foundation}}

\date{A.A.Kharkevich Institute for Information Transmission Problems \\
B.Karetnuj  per.,19, 127051, Moscow, Russia\\
 and  University of Hradec Kr\'alov\'e,
Faculty of Science, Rokitansk\'eho 62, 500~03 Hradec Kr\'alov\'e,  Czech Republic \\
\vspace{+0.50cm}
\it{\qquad \qquad \qquad \qquad  Dedicated to Jubilee of  Joseph Krasil'shchik}}


\maketitle

\abstract
   There  are  several different,  but  equivalent    definitions of  geodesics in  a Riemannian manifold,
    based  on two  characteristic  properties: geodesics  as  shortest  curves  and  geodesics  as  straightest  curves. They   are  generalized  to sub-Riemannian manifolds,  but  become  non-equivalent. We give an overview  of   different   approaches  to the definition, study and  generalisation of sub-Riemannian  geodesics and discuss   interrelations  between  different definitions.     For   Chaplygin
    transversally  homogeneous  sub-Riemannian   manifold  $Q$,
       we prove   that    straightest geodesics
      (defined  as  geodesics  of    the  Schouten partial   connection)
coincide with shortest geodesics
  (defined  as   the projection to  $Q$  of   integral  curves (with trivial initial covector)  of  the  sub-Riemannian  Hamiltonian  system).
 This gives   a Hamiltonization  of  Chaplygin  systems  in non-holonomic mechanics.\\
We  consider     a  class of homogeneous sub-Riemannian manifolds,
 where  straightest  geodesics   coincide with  shortest  geodesics,  and  give  a  description  of
   all  sub-Riemannian   symmetric  spaces    in terms of  affine  symmetric  spaces.

\tableofcontents

\section{ Introduction}

The important role of  Riemannian geometry in applications
is based on
the fact that  many  important  equations,  arising    in  mechanics, mathematical physics,  biology, economy, information  theory, image processing  etc.,   can be  reduced   to the  geodesic  equation.  Moreover,   Riemannian geometry    gives  an  effective  tool   to  investigate   the geodesic  equation  and other   equations    associated  with  the  metric (Laplace, wave,   heat  and Schr\"odinger  equations, Einstein  equation, Yang-Mills equation  etc).\\
   There  are  many   equivalent   definitions  of geodesics in  a Riemannian manifold. They  are  naturally generalised  to  sub-Riemannian   manifolds, but become  non-equivalent.\\
 H.R. Herz   remarked  that   there  are  two main approaches  to  the definition of    geodesics:
 geodesics  as {\bf shortest  curves}    based  on Maupertruis  principle of least  action (variational approach)  and
   geodesics  as  {\bf  straightest   curves}   based  on d'Alembert's  principle  of  virtual  work
    (which leads   to  a geometric description,  based on the  notion  of connection).\\
\

    We  consider   three   variational definitions  of geodesics of  a \SR mani\-fold $(Q,D,g^D)$
    (i.e.  a manifold $Q$  with  a non-holonomic  distribution $D$  and  a Riemannian metric $g^D$
     on  $D$) as (locally) shortest  curves ( Euler-Lagrange   (EL-geodesics), Pontryagin (P-geodesics)   and  Hamilton (H-geodesics)) and  three  geometric   definitions of \SR geodesics  as  straightest  curves ( d'Alembert (dA-geodesics), Schouten-Synge-Vran\-cea\-nu (S-geodesics) and
      Morimoto (M-geo\-de\-sics)) and discuss   interrelations  between  them. \\
      The definition of   M-geodesics  is based   on  E.Cartan frame bundle  definition  of Riemannian geodesics,  which is naturally generalized   to  Cartan  connections  and $G$-structures of  finite  type.
       We   give a short introduction   to this  theory in  section 4. In  section 5,  we  discuss   the  relation  between Cartan connections  and  Tanaka structures (or non-holonomic $G$-structures). They  are  defined  as  a  $G$-principal bundle  $\pi : P \to Q = P/G$ of frames  on a non-holonomic  distribution $D \subset TQ$.
      In particular, a regular \SR manifold  $(Q, D, g^D)$ (see Sect. 5.1) may be identified  with  a  Tanaka structure $\pi: P \to Q$  of admissible orthonormal frames in $D$.\\
        Using his   theory of filtered  manifold, T. Morimoto  proved  that  this Tanaka structure   admits  a unique   normal Cartan  connection, i.e.  a Cartan  connection  with  coclosed curvature.
        The Morimoto  geodesics  are  defined in  terms of  this Cartan  connection.
        We    give  a  simple  description   of  all (not necessary normal)  Cartan   connections,  associated  to  a regular \SR manifold  $(Q, D, g^D)$    in term  of  admissible  riggings  $V$ (some distribution, which is complement  to $D$)  and  define  Cartan-Morimoto (shortly, CM) geodesics  in terms  of  such Cartan  connections. CM-geodesics  are   horizontal geodesics of  some    Riemannian   connection with  torsion,  which preserves  the  distribution.
          A  necessary and  sufficient  condition  that  CM-geodesics  coincides with  S-geodesics (i.e. geodesics of the  partial \SR Schouten  connection, associated  with  a given rigging V) is  given.\\
         A. Vershik  and  L. Faddeev  \cite{V-F}   had  formulated   the  problem how to  characterize   sub-Riemannian   manifolds
  such  that   straightest  S-geodesics "coincide"(  more precisely, consistent)
  with shortest H-geodesics in  the  following  sense.\\
      An  S-geodesic  $\g(t)$ of  a   \SR  manifold $(Q,D,g^D)$   is  determined  by  the initial  velocity  $  \dot{\g}(0) \in D_q \subset  T_qQ$. The   initial    data   for an  H-geodesic    is  a pair  $(\dot{\g}(0), \l) \in   D_q \times D_q^0$  where  $D^0
        \subset  T^*Q $  is  the  codistribution (the annihilator  of   the  distribution  $D$). The  covector $\lambda$ is called  the initial codistribution covector.\\
         Taking this into account,   we  say, following  \cite{V-F},    that (straightest) S-geodesics  coincide  with (shortest) H-geodesics   if  the  class  of   S-geodesics  coincides  with   the  class of H-geodesics  with  zero  initial codistribution covector.\\
  Vershik   and Faddeev    showed  that for    generic sub-Riemannian    manifolds almost  all shortest   geodesics   are different  from  straightest  geodesics.   They   gave  the  first   example
    when shortest geodesics  coincide  with   straightest geodesics with   zero codistribution covector.\\
     In  the  second part of  the paper, we   show    that   this is  true  for any  Chaplygin  system, that is    $G$-invariant \SR   metric $(D = \mathrm{ker}(\vpi), g^D)$   on  the  total  space   of  a
 $G$-principal bundle  $\pi : Q \to  M = Q/G$  over  a Riemannian manifold $(M,g^M)$   with  a principal   connection  $\vpi : TQ \to \gg$,   where   $g^D$ is  the metric in $D$,  induced   by  the Riemannian metric  $g^M$.\\
 Any left-invariant metric on the  group $G$ defines   an  extension  of  the \SR metric $g^D$ to a Riemannian metric $g^Q$ on  $Q$.
 We  show    that H-geodesics of \SR  Chaplygin metric    are  the  horizontal   lifts  of the  projection  to  $M$  of   geodesics   of     the  Riemannian  metric   $g^Q$  and S-geodesics  are  horizontal lift  of  geodesics of  the Riemannian metric  $g^M$.
This is  a generalization  of   results  by  R. Montgomery  \cite{Mont}, 
who considered    the case  when  the  extended  metric  $g^Q$ is  defined  by  a  bi-invariant metric on  $G$.
We  give  a  simple  proof of  Wong  results  on the description  of  the  evolution  of  charge  particle  in   a  classical Yang-Mills field  in  terms of  geodesics of  the  bi-invariant    extension $g^Q$ of   the Chaplygin   sub-Riemannian metric.\\
In  the last  section,  we    describe  some  classes of invariant \SR  structures   on   homogeneous   manifolds,  where straightest  geodesics  coincide  with shortest  ones.\\
     We  give   also  a  simple  description  of  all    bracket generating symmetric \SR manifolds, introduced  by   R.S. Strichartz  \cite{Str}, and    show  that  any  flag  manifold  of  a  compact semisimple Lie  group $G$, associated  to a   gradation  of  depth $k>1$ of  the  corresponding  complex  semisimple Lie  algebra,  has  a structure of \SR  symmetric  space.\\

   Acknowledgment.  I  thank A.M. Vershik  for   his  comments  on non-holonomic   geometry  and  explanation  of  his  joint  with L.D. Faddeev  results   and  A. Spiro   for useful  discussions. 

\section{Sub-Riemannian geodesics as shortest curves}
 Here  we briefly   discuss three approaches to the definition of geodesics of sub-Riemannian manifolds :
  Euler-Lagrange variational  approach, Pontryagin  optimal control  approach  and  Hamiltonian approach.
  We  describe interrelation  between  corresponding  notions  of \SR  geodesics :  EL-geodesics, P-geodesics  and H-geodesics.

\subsection{Euler-Lagrange sub-Riemannian geodesics}

   Recall  that   a  rank-$m$ distribution $D \subset  TQ$ on  a connected  $n$-dimensional  manifold   $Q$ is  called {\bf  bracket generating }
    if    the  space $\G D$  of  sections   generates   the Lie  algebra  $\cx(M)$  of  vector  fields.\\
   According  to Rashevsky-Chow   theorem,  any  two points  on    such manifold  can be  joint  by a horizontal     (i.e. tangent  to   $D$) curve.  Then    any  Lagrangian  $L  \in  C^{\infty}(TQ) $
   defines   a nonholonomic  variational problem: \\
          Let  $C^D(q_0,q_1)$  be the  space  of  horizontal  curves,  connecting   points  $q_0$  and $q_1$. Find  a  curve  $  q(t), \,  t \in  [0,T]$  in  $C^D(q_0,q_1)$
  which delivers a minimum  or, more generally,   a critical  point,   for  the action functional
 $$A( q(t)) = \int_0^T L(q(t), \dot{q}(t)) dt , \,\,  q(t) \in C^D(q_0,q_1). $$
 The  Lagrangian  $L(q, \dot q)$ determines   a   horizontal 1-form
 $$ F_{L} = (\delta L)_i dq^i :=  (\frac{d}{dt}L_{\dot q_i}  - L_{q_i})dq^i $$
on  $TQ$, called  the  Lagrangian  force \cite{Sy}, \cite{V-F}.
 Locally    the   distribution    $D$   is  the  kernel   of   a   system  $(\o^1,\cdots, \o^k), \,  k = n-m , $ of  1-forms. The  1-form
 $$\o_{\l} = \sum \l_a \o^a =   \sum \l_a(q, \dot q)\o^a_i(q) dq^i $$
  vanishes on  $D$   for  any  vector-function
  $\l(q, \dot q) =(\l_1, \cdots, \l_m)  $ on  $TQ $.\\

 Then   critical points   $q(t)$  of the  functional $A(q(t))$ are   solution  of the  Euler-Lagrange  equations  \cite{V-G1}
 \be \label{EL equation}
  F_L\equiv (\delta L)_i\dot{q}^i(t) = \cl_{\dot{q}} \o_{\l}=  \dot{\l}_a \o^a + \l_a\dot{q} \lrcorner d\o^a
 \ee
$$    \o_{\l}(\dot{q}) =0     $$
for  unknown   curve $q(t) \in C^D(q_0,q_1)$ and   vector-function $ \l(q(t), \dot{q}(t))$.
  Here $\cl_{\dot{q}(t)}$ is  the Lie  derivative  along  the  vector  field   $\dot{q}(t)$.\\

   A   {\bf  sub-Riemannian  manifold }   $(Q,D,g^D)$ is     a manifold  with a  distribution $D$ and  a Riemannian metric  $g^D$
    on  $D$ .\\

    Let $(Q,D,g^D)$   be  a \SR manifold  with    bracket  generating   distribution $D$.
 An {\bf Euler-Lagrange   or EL  non-parametrized geodesic} ( resp., {\bf  EL-parametrized  geodesic})
   is  a  critical  point   of  the  length  functional  with $L = \sqrt{g(\dot q, \dot q)})$
 ( resp.,  the   energy  functional  with  $L = \frac12 g(\dot q, \dot q)$   )
  in  the  space $C^D(q_0,q_1)$.

\subsection{ Pontryagin  sub-Riemannian   geodesics}
 Let $(Q, D, g^D)$ be    a bracket generating \SR manifold  as  above.
  Denote  by   $(X_1, \cdots, X_m)$   a   field of orthonormal  frames in  $D$.
   Then  any    horizontal  curve $q(t) \in  C^D(q_0,q_1)$ is  a   solution  of
 the first order    ODE
   \be  \label{Control equation} \dot q(t)  = \sum_{i=1}^m  u^i(t) X_i(q(t)),\,\, \,   q(0)=q_0.
   \ee

where  
the  vector-function
$u(t) =  (u^1(t), \cdots, u^m(t))$  ( called the {\bf  control} ) consists   of  the  coordinates of the  velocity vector  field $\dot q(t)$  with   respect  to  the  frame $(X_i)$ . \\
 The  vector-function $u(t)$ is  called an  {\bf admissible  control}  if  the  solution $q(t)$ of (\ref{Control equation}) belongs  to  $C^D(q_0,q_1)$.\\
The energy  $E(u(t)) =  \frac12 \int_{0}^T  \sum u^i(t)^2 dt$ and  the  length $\ell(u) =  \int_0^T \sqrt{\sum u^i(t)^2} dt $
 of the   solution   $q(t)$ of (\ref{Control equation}) depend only on the  control $u(t)$ and  may be  considered  as  the
 functionals  (called   the  {\bf  cost functionals})   on  the  space  of  admissible  controls.


 A {\bf parametrized }  (  respectively, {\bf non-parametrized} ) Pontryagin   geodesic (  shortly, {\bf P-geodesic}) is  defined  as   the  integral  curve $q^u(t) \in  C^D(q_0,q_1 )$ of  the equation    (\ref{Control equation})   with  an admissible   control   $u(t)$,   which is  a critical point  of  the  cost
  functional  $E(u(t))$  ( respectively, $\ell(u(t))$). \\
 A P-geodesic  with  an  admissible  control $u(t)$,  which   delivers  a minimum of the cost  functional is called a {\bf  minimizer} or a {\bf minimal geodesic}.
 P-geodesics  coincide  with  EL-geodesics and  are  locally minimizers, \cite{A-S}.

\subsection{Hamiltonian   sub-Riemannian  geodesics }
Let  $(Q,D,g^D)$ be  a \SR manifold.     Denote  by  $(X_i)$ an orthonormal  frame in  $D$ and  by  $(\theta^i )$  the  dual coframe.  Then  the restriction $\xi_D $  of   a covector   $\xi \in  T^*Q$ to $D$  has coordinates  $p_i(\xi) :=  \xi(X_i)$ and  can be  written as $\xi = p_i \theta^i $.
   The  inverse   $ (g^D)^{-1}$   of   the sub-Riemannian metric $g^D$ is  a non-degenerate metric
in  the  dual  to  $D$   vector  bundle  $D^*$.
It  defines   a degenerate  symmetric  bilinear form  $g^* \in \Gamma(S^2TQ)$  in  $T^*Q$, called  the {\bf cometric}, which is  given   by
 $$g^* (\xi, \xi) = (g^D)^{-1}(\xi_D, \xi_D) = \sum_i p_i(\xi)^2,\,  \xi \in T^*Q.$$

  The  function  $ h_{g^D}(\xi, \xi)  = \frac12 g^*(  \xi, \xi)= \frac12 \sum p_i(\xi)^2$ on  $T^*Q$ is  called  the {\bf  sub-Riemannian Hamiltonian}.
   {\bf H-geodesics }  are  projection to  $Q$ of    orbits of    Hamiltonian  vector  field  $\vec{h} = \omega^{-1} dh \in  \cx(T^*Q)$    with  quadratic ( degenerate ) sub-Riemannian  Hamiltonian  $h_{g^D}(\xi,\xi) = \frac12 {g^*(\xi,\xi)}$.
            Here  $\omega  = dp_a \wedge  dq^a,\,  a = 1, \cdots, n $ is  the   standard  symplectic   form  of  $T^*Q$.

\subsection{Pontryagin Maximum Principle}
Recall  that vector fields $X = X^a \p_{q^a}$,  where   $(x^a)$ are local  coordinates in  $Q$ bijectively  correspond  to
  fiberwise linear   functions
   $$  p_X : T^*Q \to \bR,\,   \xi = p_a \p_{q^a} \mapsto  p(X) = X^a p_a$$
on  the  cotangent  bundle  $T^*Q$. The  function $p_X$  is  the  Hamiltonian  of the Hamiltonian  vector  field
 $$\vec{p}_X =  X^a \p_{q^a}  - p_a \p_{q^b}X^a\p_{p_b},$$
  which is  the  complete lift of  $X$ to  $T^*Q$.  The   map
   $$\cx(Q) \ni  X \to  \vec{p}_X \in \cx(T^*Q) $$
  is   an  isomorphism of    the Lie  algebra $\cx(Q)$ of vector  field  onto  the  Lie  algebra   $\cx(T^*Q)^1 $ of fiberwise linear  vector   fields  on   $T^*Q$. \\

\bt  (Pontryagin  Maximum Principle)
Let  $q(t)  \in C^D(q_0,q_1)$  be  a minimal  P-geodesic  on  a sub-Riemannian manifold
$(Q,D = \mathrm{span}(X_i),g^D)$ with  natural parametrization   (s.t.    $|\dot{q}(t)| =  const $), which  corresponds to   a  control $u(t) = (u^i(t))$ :
\be  \label{Equation for curve}
 \dot q(t) = u^i(t)X_i(q(t)).  \ee
Denote  by $\varphi_t$  the  (local) flow, generated  by  the  non-autonomous vector  field \\
$X^u =u^i(t)X_i $.
Then   for  some covector $\xi_0 \in  T_{q_0}^*Q$ the   curve
$$\xi(t) := \varphi_{-t}^* \xi_0 := \xi_0 \circ \varphi_{-t_*} \in  T^*_{q(t)}Q  $$
 satisfies the  equation
\be  \label{Equation for  extremal}
  \dot{\xi}(t) = u^i(t)\vec{p}_i(\xi(t))
    \ee
  where  $p_i := p_{X_i}$   and one of  the following  conditions  holds
$$
\begin{array}{cccc}
               u^i(t)& \equiv & <\xi(t), X_i(q(t))>     \qquad  & (N) \\
                    0 &\equiv &  <\xi(t), X_i(q(t))> .      \qquad  & (A)
\end{array}
$$
\et

Here  the  bracket  $<\xi, X>$ denotes  the pairing  between   covectors  and  vectors.\\

   An  extremal    curve   $\xi(t) \subset  T^*Q$,  which  satisfies  $(N)$ (  resp., $(A)$), is  called a {\bf normal} (  resp., an {\bf abnormal}) {\bf extremal},  and  its  projection $q(t) \subset Q$  is  called  a {\bf normal } (resp., an {\bf abnormal})  {\bf  P-geodesic}.
Note  that  abnormal extremals   are  curves   in  the  codistribution   $D^0$,  considered  as  a  submanifold  of  $T^*Q$.

    \subsubsection{  Normal  P-geodesics  as  H-geodesics}
Pontryagin  theorem  shows  that     normal geodesics  are  H-geodesics. More precisely,  we have
\bc Let $D$  be a rank-$m$ bracket generating  distribution  with  a  sub-Riemannian metric $g^{D}$.
A normal  extremal $ \xi(t) \subset   T^*Q$    for   $(Q,D,g^D)$ is an integral  curve of  the  Hamiltonian equation  on   $T^*Q$ with    the  sub-Riemannian  Hamiltonian
 $$h_{g^D}(\xi) = \frac12 g^*(\xi,\xi)  = \frac12 \sum_{i=1}^m p_i(\xi)^2 ,$$
where  $\xi = p_i(\xi) \theta^i, \,\, p_i(\xi) = \xi(X_i) $.
\ec

\pf  In  the  case  of normal geodesic, the  equation
(\ref{Equation for  extremal})  take  the  form
$$  \dot \xi(t)= \sum p_i(\xi(t)) \vec{p}_i (\xi(t))=\frac12 \o^{-1}d  \sum p^2_i(\xi(t))=\vec{h}_{g^D}(\xi(t)).
$$
\qed \\

Since  the  Hamiltonian vector  field   $\vec{h}_{g^D}   $   preserves    the  Hamiltonian  ${h}_{g^D}$,
   a  normal  extremal $\xi(t)$  belongs   to   a level  set   $L_c = \{h =c\} \subset  T^*Q$ of  the Hamiltonian
    $ h= {h}_{g^D}$.\\

 A curve $\xi(t)\subset  L $  on  a  submanifold  $L \subset T^*Q$ is  called  a {\bf  characteristic} if  its  velocity  $\dot{\xi}(t) $ belongs  to   the  kernel  $\ker (\o|{L})$ of  the   restriction  of  the  symplectic  form  to  $L$.

 \bc Assume  that an extremal  $\xi(t) \subset   L_c$ belongs  to  a  regular  level  set of     the  Hamiltonian   $h =h_{g^Q}$, i.e.  $L_c$ is  a  smooth hypersurface. Then  $\ker \o|_{L_c}$   is the 1-dimensional  distribution  generated   by   $\vec h$. In  particular,  the  extremals  $\xi(t)$ are  the  characteristic  curves   of  $L_c$.
\ec

\pf The  tangent  space  of  the  level  set  $L_c$ is  described  as   follows
 $$T_{\xi}L_c = \{ w \in  T_{\xi}(T^*Q), \,\, 0=  < dh_{\xi},w> = <\o(\vec h_{\xi}), w > = \o(\vec h_{\xi}, w) \}. $$
  This  shows  that $T_{\xi}L_c $  consists of  all $\o$-orthogonal  to   $\vec{h}_{\xi}$  vectors.
Since  the    $\o|{T_{\xi}L_c}$  has  1-dimensional  kernel, it is  generated by   $\vec{h}_{\xi}$.
\qed \\

\subsubsection{ \,  Abnormal   P-geodesics}
Now we  shortly   discuss  main   properties   of  abnormal  geodesics, following  R. Montgometry.
 Denote  by
 $C^D(q_0) = \{ \g : [0,T] \to   Q, \, \g(0) =q_0,  \dot{\g} \in  D \} $
  the  space of   horizontal curves,  starting  from  $q_0$ ,  where $D \subset  TQ$ is  a  bracket generating distribution.  A   curve $\g(t) \in  C^D(q_0)$ is  called  {\bf  singular }( resp.,{\bf regular}) if  the  end-point map
  $$\e : C^D(q_0) \to  Q, \,   \g(t) \mapsto   \g(T) $$
  is  {  singular} (  resp., {regular}).\\
The  following   theorem  by L.S. Pontryagin, L. Hsu  and  R. Montgomery  shows  that  abnormal geodesics  coincide with  singular  curves    and  they    are projection on  $Q$ of characteristic  curves   of    the  codistribution  $D^0$,  considered  as  a  submanifold  of the  symplectic manifold  $(T^*Q, \o)$.
 \bt ( see \cite{Mont}, \cite{Mont1} )
  i) Abnormal  geodesics  of  any  sub-Riemannian  metric  $g^D$  on  $D$   are   exactly    singular horizontal  curves  in $Q$.\\
     ii)  A horizontal  curve   $\g \subset  Q$  is  singular if  and only if   it is  a projection  to  $Q$  of  a  characteristic  curve   of the  submanifold   $ \tilde{D^0}:=   D^0 \setminus \{ \text{zero  section}   \}\subset  T^*Q$.\\
    \et

Now    we give  a  description  of  characteristic  curves in  $ \tilde{D^0}$, following  \cite{Mont}.  To  simplify notation,  we   will  denote
 $\tilde{D}^0$ by  $D^0$.\\

 Denote  by   $   \tau :  T^*Q \to Q,\,   \tau_*: T (T^*Q) \to  T^*Q$  the  natural projections.

 We   fix  a  complementary  to  $D$  distribution $V$  such  that $TQ = D \oplus V$.
  Let $(X_i),\,  i=1, \cdots, m$   be  a local frame in $D$, $(Y_{\a}),\,  \a = 1, \cdots, n-m$  a local  frame   in  $V$   and  denote  by     $(\theta^i_Q ,  \eta^{\a}_Q)    $   the    dual  coframe, such  that
  $$   \theta^i_Q(X_j) = \delta^i_j, \eta^{\a}_Q(Y_{\beta})= \delta^{\a}_{\beta}, \theta^i_Q(Y_{\a}) = \eta^{\a}_Q(X_i) =0.$$

The  Liouville tautological   1-form  $\theta_{\xi} =\tau^*  \xi = \xi \circ \tau_*  $  in $T^*Q = D^* \oplus D^0$
at a point
 $\xi = h_i \theta^i_Q + k_{\a} \eta^{\a}_Q   \in   T^*Q$ can be written  as
$$  \theta_{\xi}  = h_i \theta^i + k_{\a} \eta^{\a} \in  (\tau_*)^* D^*_{\xi}  \oplus (\tau_*)^* D^0_{\xi} \subset T^*( T^*Q),$$
where   $\theta^i =  \theta^i_Q \circ \tau_*, \,  \eta^{\a}  = \eta^{\a}_Q \circ \tau_*$  are   the pull back of the 1-forms $\theta^i_Q, \, \eta^{\a}_Q$   to    $T^*Q$.\\
We  will  consider  $h_i, \, k_{\a}$ as  fiberwise  coordinates  in  the  bundle  $T^*Q$  and  in  the  bundle
 $\tau^*(T^*Q) \subset T^*(T^*Q)$ of  horizontal  1-forms on $T^*Q$.\\

The   restrictions $\theta^0, \, \o^0 $   of  the  Liuville  form   $\theta$  and  the standard symplectic  form   $\o = -d \theta $  on  $T^*Q$ to  the  submanifold  $D^0 \subset T^*Q$   are given  by
$$\theta^0_{\xi} = \xi|_{D}  =k_{\a} \eta^{\a},\,\,\,
  -\o^0 = d \theta^0 = dk_{\a} \wedge  \eta^{\a}  + k_{\a} d \eta^{\a}.$$

Denote  by
 $$  Ch(D^0)  := \ker \o^0 =  T(D^0)^{\perp} \cap T(D^0)= \{v \in T_{\eta}D^0,\,  \o(v, T_{\eta}D^0) =0\} $$
  the characteristic  submanifold   of $T(D^0)$,  where   the  vector  bundle
   $T(D^0)^{\perp}  \subset T(T^*Q)|_{D^0}$  is  the   $\o$-orthogonal  complement  to  the  tangent  bundle   $T(D^0) $   of  $D^0  $.\\
  The fiber   $Ch_{\eta}(D^0)= \ker \o^0_{\eta} \subset T_{\eta}(Q^0)$ over  a  point  $\eta \in Q^0$ is a  vector  space,  but since  the rank of  $\o^0_{\eta}$ may vary, the natural  projection  $Ch(D^0) \to   D^0$ is not  a vector  bundle, in general.\\

 By  definition,{\bf characteristic  curves}  are
    curves   $\eta(t) \subset D^0$, tangent  to  the characteristic manifold $Ch(D^0) \subset  T(D^0)$.


  \bl  \label{isomorphismtau} The  vector  bundle $T(D^0)^{\perp}  =  \mathrm{span} \{ \vec{h}_i , \, i=1, \cdots, m \} $,  and   the  projection $\tau_* : T(T^*Q) \to  TQ)$
induces  for   any  $\eta \in  D^0$ the    isomorphism
$$   \tau_* : T_{\eta}(D^0)^{\perp} \to D_q, \,  q= \tau(\eta), $$
$$ u^i \vec{h}_i \mapsto u^i X_{i}|_q  .                          $$
\el
\pf The  submanifold $D^0 = \{ \eta = k_{\a}\eta^{\a}   \}$ is  defined  by the equations
 $$     h_i =0, \, i= 1, \cdots, m.,  $$
Hence,
 $$T_{\eta}D^0 = \{   v \in T_{\eta}(T^*Q) , 0 = <dh_i, v> = \o(\o^{-1}dh_i, v) =  \o(\vec{h}_i, v)     \}.$$
\qed

Since   $D^0 = \{  \eta = k_{\a} \eta^{\a}\}$,
 $k_{\a}$
 are  fiberwise  coordinate  of  the bundle  $D^0 \to  Q$.  We identify    $X_i, Y_{\a}$ with "horizontal" vector fields in
 $T^*Q$,  which  annihilate   the    fiberwise  coordinates   $h_i, k_{\a}$.   Then
   $\p_{k_{\a}}, X_i, Y_{\a}$   form     a  frame  in the  tangent  bundle
    $T(D^0)$.    The  tangent  vector   to   a  curve
    $\h(t)= k_{\a}(t) \eta^{\a}(t)  \subset D^0  $  with  projection  $ \gamma(t)= \tau  \eta(t)$ can be  written  as
   \be \label{charcurve}
    \dot{\h}(t)  =  \dot{k}_{\a}\p_{k_{\a}}  + \dot{\g}^i X_i(\g(t))+ \dot{\g}^{\a} Y_{\a}(\g(t)).
    \ee

   We   need  an  explicit    description  of  the  restriction     $-\o^0 = dk_{\a}\wedge \eta^{\a} +  d \eta^{\a}$  to $D^0$   of  the  symplectic   2-form   $\o = - d \theta$.
   We   may  write   the   2-form $ d \eta^{\a}$   as
  $$   d \eta^{\a} = -c^{\a}_{ij} \theta^{i}\wedge  \theta^{j}-c^{\a}_{i \beta} \theta^{i}\wedge\eta^{\beta}  -c^{\a}_{\beta \delta} \eta^{\beta}\wedge  \eta^{\delta}
    $$
  where  $c^{\a}_{ij} = \eta^{\a} ([X_i, X_j]), \,c^{\a}_{\beta i} = \eta^{\a} ([Y_{\beta}, X_i]), \,
   c^{\a}_{\beta  \delta} = \eta^{\a} (Y_{\beta}, Y_{\delta}).$\\
Then
$$  \o^0 = -d k_{\a} \wedge \eta^{\a} + k_{\a}( c^{\a}_{ij} \theta^i \wedge \theta^j + c^{\a}_{i \beta} \theta^i \wedge \eta^{\beta}  +  c^{\a}_{ \beta \delta} \eta^{\beta} \wedge \eta^{\delta} ).                 $$
Now  we  are ready  to   write  down  the  necessary and  sufficient  condition that a  tangent vector $\dot{\eta}(t) \in T_{\eta(t)}D^0$  of  a  curve  $\eta(t) \subset D^0$  belongs  to $Ch_{\eta(t)}(D^0)$.\\
 Calculating   the  contraction  ${\dot{\eta}}\lrcorner \o^0  $,  we   get
 $$   {\dot{\h}}\lrcorner \o^0 = \dot{\g}^{\a} dk_{\a}   + \k_{\a} (c^{\a}_{ij}\dot{\g}^j -c^{\a}_{i \beta} \dot{\g}^{\beta})\theta^i + ({k}_{\a}c^{\a}_{i \beta} \dot{\g}^i + {k}_{\a}c^{\a}_{ \beta \delta} \dot{\g}^{\delta}-\dot{k}_{\beta} )\eta^{\beta}.                              $$
  In particular, a  curve   $\eta(t ) =  k_{\a} \eta^{\a} \subset   D^0  $ is  a  characteristic curve  if and only if  its  velocity  vector (\ref{charcurve})  satisfies  the  following   equations
$$
 \begin{array}{llll}
i)&   & \dot {\g}^{\a} & =0, \\
 ii)& &  \dot{k}_{\beta} - k_{\a}c^{\a}_{i \beta}\dot{\g}^{i}& =0,\\
iii)&  &    k_{\a}c^{\a}_{ij} \dot{\g}^i& =0.
 \end{array}
 $$

For  $q \in Q$, denote  by   $\Lambda^2 D^*_q $  the  space of  2-forms   in $D_q$  and  by $D^0_q$   the  fiber  of  the  bundle   $\tau :  D^0 \to Q$.  There is   a natural linear  map
  $$\bar{d} : D_q  \to \Lambda^2 D^*_q$$
    $$  \eta  \mapsto   d \tilde \eta|_{\Lambda^2D_x},  $$
   where    $\tilde \eta$ is  an    extension    of       $\eta$   to a  local  1-form.   If  $\tilde{ X}, \tilde{X}'$ are    extensions   of   vectors  $X,X' \in  D_q$ to  local  sections of $D$,  then
   $$  \bar d\eta(X,X')   = - \eta ([\tilde X, \tilde X']).$$
     This  shows  that  the   map $\bar d$  does  not  depend  on  extensions  $\tilde{\xi}, \tilde{X}, \tilde{X}'$.\\
 We   set
  $$K_{\eta} = \ker \bar{d}\eta \subset D_{\tau(\eta)}.$$
  \bp  The  projection  $\tau_* : T_{\eta}(D^0) \to  D_{\tau(\eta)}  $  induces  an isomorphism
   $$   \tau_* : Ch_{\eta}(D^0) = \ker \o^0_{\eta}  \to K_{\eta} .             $$
  \ep

  \pf Lemma \ref{isomorphismtau} shows  that   $\tau_* : T_{\eta}D^0 \to D_{\tau(\eta)}$ is  an isomorphism.
  The   conditions i), iii)   may  be  rewritten  as
  $$ \dot{\g}=  \dot{\g}^i(t)X_i    \in K_{\eta(t)}=  K_{k_{\a}\eta^{\a}}.$$

     Any  characteristic  vector   $\dot{\eta} \in Ch_{\eta}(D^0) = \ker\o^0_{\eta}$
  can be  written now  as
  $$   \dot{\eta} = \dot{k}_{\a}\p_{k_{\a}} + \dot{\g} = k_{\beta}c^{\beta}_{\a i}\dot{\g}^i \p_{k_{\a}} + \dot{\g}^i X_i                   $$
   and it is   completely  determined  by   the  point  $\eta = k_{\a}\eta^{\a}  \in D_q^0$  and the tangent  vector
    $\dot{\g} \in K_{\eta} \subset D_{\tau(\eta)} $.
   \qed

   As  a  corollary,  we  get  the  following  characterization  of characteristic  curves  and   abnormal  geodesics.
   \bt
   i) A  curve $\eta(t)= k_{\a}(t) \eta^{\a}( \gamma(t)) \subset D^0$ with  the projection  $\gamma(t)=\tau (\eta(t)) $
   is  a characteristic  and  then  $\g(t)$ is  an  abnormal geodesic if  and only if
   the velocity  vector  field  has  the  form
    $$  \dot{\eta}(t) = k_{\beta}c^{\beta}_{\a i}\dot{\g}^i \p_{k_{\a}} +\dot{\g}^i X_i  $$
   such  that  $\dot{\g}(t) =  \dot{\g}^i X_i \in K_{\eta(t)}$.\\
     ii)  A horizontal  curve  $\g(t) \subset Q$  with  velocity  vector  field  $\dot \g(t) = \dot{\g}^i(t) X_i(\g(t))$
    is  an abnormal geodesic  if  and only if  it  can be  lifted  to a characteristic  curve  $\eta(t)  \subset  D^0$
    such  that  $\dot{\g}(t) \in K_{\eta(t)}$.\\
   \et

\section{ Sub-Riemannian geodesics  as   straightest  curves}

\subsection{  d'Alembert's sub-Riemannian  geodesics} 
Let  $(Q,D,g^D)$ be a \SR manifold. To define d'Alembert's (shortly, dA) geodesics, we  extend the  sub-Riemannian  metric $g^D$ to a Riemannian  metric $g^Q$.\\
 The  d'Alembert's principle  of  virtual displacements  for  a mechanical  system  may  be  formulated  as  follows, see \cite{V-F}.\\
1) The   evolution of  a mechanical  system  with  a (smooth)  configuration  space  $Q$ is  described  by projection  to   $Q$  of  integral  curves   of  a  special  vector  field  $X  \in \cx(TQ)$    
(the evolution  field). A field $X$ is called {\bf  special}  if it corresponds  to   a second order equation,  that is $\pi_* X_{(q, \dot q)}  = \dot q$   where    $\pi:TQ \to Q$ is  the projection.\\
2) The  vector  field   $X$ is  determined  by  the Lagrangian  force,
  defined  as  the  horizontal  1-form  $ F_L:=(\delta L(q,\dot q))_i  dq^i $ on  $TQ$,  associated  with  the  Lagrangian  $L(q,\dot q)$,  and
   external  forces.\\
 3)  d'Alembert's Principle  states  that  the special vector  field  $X$, which  describes  the  real  dynamics of  a  mechanical   system,
  is  determined  by  the  condition  that  the  Lagrangian  force is equal  to   the  external  force.\\
  Assume   that\\
   i) the Lagrangian $L(q, \dot q)$ of  the  system  with    a configuration  space   $Q$ is  quadratic in  velocities $\dot q$ and  positively  defined  (  that   is can be  written  as $L = \frac12 g(\dot q, \dot q)$,  where  $g$ is  a Riemannian metric in  $Q$ )     and  that \\
  ii)    the  only   external  force   is  the    reaction   of a non-holonomic  constraint,  defined   by a rank-$m$  distribution   $D  =   \ker \eta^1 \cap \cdots \cap \ker \eta^k$,  where   $\eta^{\a} = \eta^{\a}_i dq^i,\,  \a= 1, \cdots, k=n-m$ is  a  coframe  of   the  codistribution  $D^0$. The   reaction of  the  constraint  is  the horizontal  1-form    $\phi_{\l} = \l_{\a} (q, \dot q)\eta^{\a}  \in \Omega^1(TQ)  $    defined  by   the condition    that the   equation
  $F_L - \phi_{\l} =0$  corresponds  to a vector  field  $X \in  \cx(TQ)$
       tangent  to the  distribution  $D \subset  TQ$.\\
        In  coordinates,  this  equation   take the form \cite{V-G1}
   $$   F_{L}\equiv  (\frac{d}{dt}L_{\dot q_i}- L_{q_i})dq^i = \l_{\a}(q(t),\dot{q}(t))\eta^{\a}_i dq^i $$
or
   $$     \frac{d}{dt}L_{\dot q_i}- L_{q_i} \equiv 0 \quad (\mathrm{mod}  D^0).$$
  The  projection  to $Q$ of  integral  curves  of  this  equation  is  called {\bf dA-geodesic } of the \SR metric   $(D, g^D)$,  associated  with   an  extension  of  $g^D$  to  a Riemannian metric  $g$  on $Q$.
    In  general,    the equation   of  dA-geodesics  is  neither  Lagrangian  nor  Hamiltonian.\\

\subsection{Schouten-Synge-Vranceanu sub-Riemannian geodesics}

     Recall that  Levi-Civita
         associated  to  a  Riemannian manifold $(Q,g) $      the  canonical   torsion  free connection $\n^g$,  which preserves the  metric  (called   the  Levi-Civita  connection).   According   to Levi-Civita, a geodesic  is  defined  as an autoparallel curve  $q(t)$, such  that  the velocity vector  field  $\dot{q}(t)$  is parallel  along $q(t)$, i.e.   satisfies  the geodesic  equation

  $$ \n^g_{\dot \g}\dot \g \equiv \ddot q^i(t) + \G^i_{jk}(q^j(t)) \dot q^j(t) \dot q^k(t)  =0 $$
where
 $\G^i_{jk}$  are the Christoffel  symbols of  the metric  $ g =  g_{ij}(q) dq^i dq^j$.
 The  extension  of  this  definition   to   \SR manifolds  had  been  proposed independently  by   J.A. Schouten,   J.L.Synge and  G. Vranceanu,  see   \cite{B}.\\
\subsubsection{Schouten partial  connection of a sub-Riemannian manifold}
   Let    $D \subset  TQ$ be  a  distribution.  A partial  $D$-connection  in  $D$ is an  $\bR$-bilinear  map
  $$ \n^D :  \G D \times \G D \to  \G  D, \, (X,Y) \mapsto \n^D_X Y $$
   which  is
  $ C^{\infty}(Q)$  linear  in  $X$ and  satisfies   the  Leibnitz  rule  in $Y$:
  $$\n^D_X (f Y) = f \n^D_X Y    + (X\cdot  f)Y,\,\,    f \in C^{\infty}(Q). $$

    Let $e_i, \,  i = 1, \cdots, m$  be  a  frame  of  $D$   defined in  a neighborhood    of   a  horizontal  curve $q(t)$. \\
   The Christoffel  symbols of  a  partial  connection $\n^D$  are  the local  functions   $\G^i_{jk}(q) $ on $Q$ defined  by
  $$  \nabla_{e_j}e_k  =   \G^i_{jk}(q) e_i. $$
 The   value  of  the   functions  $\G^i_{jk}(t) := \G^i_{jk}(q(t))$ on a horizontal   curve $q(t)$ depends
 only on  the   frame   $e_i(t):=e_i( q(t))$  along  the  curve  $q(t)$.
 Due  to  this,
 the   partial  connection  defines a parallel transport   of a  vector   $Y_0 \in D_{q_0}$  along a  horizontal  curve
$  q_t  $   as    the   solution $Y(t)=  Y^c(t)e_c(t) \in D_{q_t}$   of  the  equation
$$ 0=  \n_{\dot{q}_t} Y(t) =\n_{\dot{q}_t}(Y^i(t)e_i(t))  = [\dot{Y}^i(t) + \Gamma^i_{jk}(t) q^j_t Y^k(t)]  e_i(t).$$   
I.A. Schouten   showed   that  a  complementary  to  $D$  distribution  $V$  on  a \SR  manifold    $(Q,D,g^D)$  ( called   a {\bf rigging})  defines  a partial connection  $\n^S$ in   $D$    which   preserves   the  metric $g^D$    and  has  zero  torsion $T$.  The torsion tensor   is  defined  by
 $$T(X,Y) = \n^S_X Y - \n^S_YX - [X,Y]_D, \,  X,Y \in  \G D ,$$
 where  $X_D$  is  the  horizontal part  of  the   vector
 $$X =  X_D + X_V \in  T_q Q = D_q \oplus V_q.$$

In  coordinate-free  way, the Schouten  partial  connection  of  $(Q,D,g^D)$  associated  to a rigging  $V$  is defined  by  the Koszul  formula
 $$
  \begin{array}{ll}
   2g(\n^S_X Y, Z) =& X \cdot g(Y,Z) +  Y \cdot g(X,Z)- Z\cdot g(X,Y)+\\
                   &  g([X,Y]_D,Z)  -  g(Y, [X,Z]_D)  -  g(X, [Y,Z]_D),\\
                   & X,Y,Z  \in  \G(D).
    \end{array}
   $$
   Schouten defined  the  curvature tensor  $R  \in \mathfrak{so}(D)\otimes \Lambda^2 T^*M $  of  the Schouten  connection    by
$$        R(X,Y)Z =  [\n_X, \n_Y]Z - [[X,Y]_V,Z ]_D, \,  X,Y,Z  \in  \G D. $$

V.V. Wagner generalized this  notion  and  defined Wagner curvature  tensor,  such  that  the   vanishing of  the  Wagner  tensor  is equivalent   to  the flatness  of  the Schouten  connection ( that is the property  that  the    associated  parallel  transport    does not   depend on  the  path, connecting two points),  see \cite{B}, \cite{D-G} for  a modern  exposition  and generalization  of  this  theory.\\

\subsubsection{Sub-Riemannian   S-geodesics  and non-holonomic mechanics}
{\bf   Schouten-Synge-Vranceanu   geodesics ( S-geodesics)} of     a  sub-Riemannian manifold   $(Q,D,g^D)$, associated  to a rigging $V$,
     are defined as   horizontal    curves $\g(t)$ with   parallel (w.r.t.   Schouten  connection) tangent  vector field  $\dot{\g}(t)$, i.e.  solutions   of  the  equation
   $  \n^S_{\dot{\g}}\dot{\g} =0.$\\
Assume  that   the   \SR   metric  $g^D$  is  extended  to a Riemannian metric
$g$ on  $Q$.  Denote  by    $ V= D^{\perp}$  the  $g$-orthogonal  complement  to  $D$.
 Then     the Levi-Civita  connection  $\n^g$
 induces   a   connection  $\n^D $   in  $D$  given  by
$$           \n^D_X Y = \mathrm{pr}_D \nabla^g_X Y = (\n^g_XY)_D , \,  X  \in TQ, Y \in \G D. $$
  where  $\mathrm{pr}_D : TQ = D \oplus D^{\perp}\to D$ is  the  natural projection.

The  connection $\n^D$  is  an  extension  of  the  partial  Schouten  connection  $\n^S$  associated to  the  rigging $D^{\perp}$.

\bt (   Vershik-Faddeev\cite{V-F}, \cite{V-F1})  Let $(Q,D,  g^D)$ be  a \SR  manifold,  $g$    an   extension of  $g^D$  to  a metric  in  $Q$
   and   $V  =  D^{\perp}$  the  orthogonal  complement  to  $D$.  Then   S-geodesics  coincide  with  dA-geodesics   and
    they  describe  evolution  of the  free  mechanical  system  with  kinetic    energy
   $  g$  in  configuration  space  $Q$  with  nonholonomic linear constraint  $D$.
\et

\section{ Cartan   frame   bundle   definition  of  geodesics }

An  important  frame bundle  definition of  Riemannian   geodesics  
  had  been  proposed  by  E. Cartan. It is  naturally generalized  to a  wide   class   of  geometric  structures.
Below  we    consider  Cartan approach  to definition  of geodesics  for   Cartan  connections     and $G$-structures of  finite  type.
   This      will be used  for  definition   of  Morimoto geodesics  on  a  regular  \SR  manifold.
\subsection{    Cartan  definition  of  Riemannian geodesics}

A Riemannian  metric $g$ on  a manifold $Q$    can be  considered  as  a $G= O_n$-structure, i.e.  a principal $G$-subbundle
   $\pi : P \to Q  =  P/G$     of the  bundle  of   orthonormal frames  (i.e. isometries
 $f: \bR^n =V \to T_xQ  $  )    with    the  tautological  soldering  form
$$\theta : TP \to V, \, \, \,\, \theta_f(X) := f^{-1}(\pi_* X).      $$

    The   total  space  $P$    of an $O_n$-structure  admits  a canonical $O_n$-equivariant   absolute parallelism (Cartan  connection)
 $$\kappa = \theta  \oplus  \omega :  TP  \to   V \oplus \mathfrak{so}(V),$$
  which is    an  extension  of  the  vertical parallelism $i_p : T^v_p P \simeq \gso_n, \,  \forall  p \in P$
  (defined  by  the  free  action  of  $O_n$   on $P$).
    Here   $T^vP \subset  TP$  is  the  vertical  subbundle and    $\omega : TP \to \gso_n$ is  the   connection  form of  the  Levi-Civita  connection.\\
   {\bf Cartan   geodesics (C-geodesics)}    are  defined    as   the   projection   to   $Q$  of    constant horizontal   vector  fields   $X \in \kappa^{-1}(V)  \subset \cx(P)$,  see   \cite{K-N}.

\subsection{  Normal Cartan  connection  and  C-geodesics \label{Cartan connections}}

We  recall  the  definitions  of   a Cartan connection  and   associated  C-geodesics.\\
Let   $M_0 = L/G$ be  an $n$-dimensional  homogeneous  manifold.  \\
A {\bf Cartan  connection  of  type     $M_0 = L/G$    }   on an $n$-dimensional manifold   $Q$    is  a principal $G$-bundle   $ \pi : P \to Q =P/G $ together  with an  $\gl$-valued  $G$-equivariant
( s.t.  $ r_g^* \kappa   =  \Ad^{-1}_g \circ \kappa, \,  g \in G$ )  kernel-free  1-form
 $\kappa :  TP  \to \gl$
   which  extends  the vertical parallelism  ${T_p^vP} \simeq \gg $.\\
   The   form $\kappa$ defines  an  absolute parallelism   $\kappa_p: T_pP \simeq \gl$.  Hence,   tensor  fields on $P$ may be  identified  with  tensor-valued  functions.\\
In  particular,  the  horizontal  {\bf curvature 2-form }
$ \Omega:= d \kappa  + \frac12 [\kappa, \kappa]$
 on  $P$       can be  identified   with   a function
 $ K : P \to  C^2(\gn, \gl) :=  \Lambda^2 \gn^*  \ot \gl$
     where  $\gn  = \gl/\gg$.\\

 One of  the  most powerful method for  studying     different  (holonomic  and non-holonomic) geometric  structures  and  for constructing  their  invariants  is based  on   construction of  the associated  canonically defined   Cartan  connection.  In many cases, it is not  difficult  to   associate  to   a  given  structure   a   family  of Cartan  connections. Then
the problem comes down to finding suitable normalization conditions    which  uniquely specify a  Cartan  connection (called the
{\bf normal Cartan  connection}) of this  family.    The  standard  way  is  to  impose  some  normalization   conditions
 on  the  curvature  function $K$,  for  example,  the  condition  that    the  curvature  tensor $K_p,\, \forall p \in  P$ is  coclosed.\\
    We    explain   this  condition   in  the case  when the  homogeneous manifold  $M_0 =L/G$  satisfies the  following property :\\

    $(*) $  The  Lie  algebra  $\gl$ admits an  $\ad_{\gg}$-invariant  metric $g$   and  the $g$-orthogonal  complement
      $\gm$ to $\gg$   in $\gl$ is a  subalgebra.\\
 This  is  sufficient  to   define  Morimoto  geodesics  for  regular  \SR  structures.\\
     Let
  $$   \p : C^1(\gm, \gl) = \mathrm{Hom}(\gm, \gm+\gg) \to C^2(\gm, \gl)            $$
be   the  differential  of the complex  of exterior  forms on  the Lie  algebra  $\gm$  with values  in  $\gm$-module  $\g$.
  Denote  by
       $$  \p^* : C^2(\gm, \gl) \to C^1(\gm, \gl)    $$
     the   dual  codifferential, defined  by  means of    the induced  metrics on $C^j(\gm, \gl)$.\\
  A {\bf  Cartan  connection is called normal }  if   the  curvature  $K$  is  coclosed, i.e.  $\p^* K =0$.\\
   \bt (see \cite{M}, \cite{A-D}, \cite{C})  Let $(\pi: P \to Q, \, \kappa : TP \to \gl = \gg + \gm)$ be  a Cartan  connection,  which  satisfies  the  condition  (*). Then the bundle  $\pi$ admits  a unique  normal  Cartan  connection  $\kappa_0$.
    \et
       More  general    sufficient  conditions  for the existence of a  unique normal Cartan  connection   are  given in  \cite{M1}, \cite{C-S},
        \cite{A-D}, \cite{C}, \cite{C-D-K}.

Let   $(\pi: P \to \gl, \,\, \kappa :TP \to  \gl)$ be  a Cartan  connection  and $\gm$     a  fixed complementary to $\gg$
subspace  of  $\gl$.  Then  $\kappa^{-1}(\gm) \subset   TP$ is  a  complementary  to $T^vP$  distribution, called  the {\bf
horizontal  distribution} and any  vector $v \in \gm$ defines  a horizontal vector  field  $X^v = \kappa^{-1}(v)$, called  the {\bf constant  horizontal vector  field associated  to  $v$}.\\

             Like in Riemannian  case,  {\bf  C-geodesics  of  a Cartan  connection}  are defined  as the projection  to  $Q$ of  integral  curves  of   constant   horizontal vector  fields.

 Assume  that  the  homogeneous  manifold  $M_0 =L/G$ is {\bf  reductive}, i.e.   there  is   a  reductive  decomposition   $\gl = \gg \oplus V$,
   where  $V$ is  an  $\Ad_G$-invariant   complement  to  $\gg$.
   Let   $(\pi, \kappa)$  be  Cartan  connection  of   a   reductive  type $M_0 = L/G$.
    Denote  by
    $$\theta := \mathrm{pr}_V \circ \kappa :  TP \to V  \,\, ( \text{ resp.},\, \o := \mathrm{pr}_{\gg}\circ \kappa :  TP \to \gg   )$$
    the  horizontal part ( resp.,  the vertical part) of  the 1-form  $\kappa$.
     Then $\theta$  is a  soldering  form,  which  turns  $\pi$ into  a $G$-structure,
      and   $\o$ is  a  connection  form,  which  defines  a  principal  connection  in  $\pi$. The form $\o$  defines
       a linear    connection   $\n$   in  the  tangent bundle   $TQ = P \times_G V$,  see \cite{K-N}  and
         C-geodesics of the Cartan  connection  coincide  with   geodesics of  $\n$.

\subsection{C-geodesics  for  $G$-structures}

\subsubsection{ $G$-structures  and  their torsion  function}
 We   recall  the definition  of $G$-structure   and   its  torsion  function.\\

Let $G \subset GL(V), \, V = \bR^n$ be  a linear Lie  group. A {\bf $G$-structure } on  an  $n$-dimensional manifold $Q$  is
  a  $G$-principal bundle $\pi  : P \to  Q =  P/G$  with  a  {\bf soldering  1-form}    $   \theta  :  TP \to  V $
i.e. a  strictly horizontal   ( $\ker \ \theta =T^vP $) $G$-equivariant   1-form.\\
Such   1-form    allows  to  identify  the $G$-principal  bundle  with   a  $G$-principal bundle of   frames  on  $TQ$. \\
  Indeed, the  soldering  form   at  a point $p \in P$  defines a  coframe, i.e.  an isomorphism
  $$ \theta_p : T_{\pi(p)}Q  \to  V. $$
 We  denote  by
  $$\hat{p} = \theta_p^{-1} :  V \to T_{\pi(p)}Q $$
    the  dual  frame. This  allows to identify  the  bundle   $\pi$ with a $G$-principal bundle of  frames.\\
       Denote  by   $ j^1(\pi):  J^1 \to P$   the  bundle  of 1-jets of local sections  $H=H_p = j^1_p(s)$,
        that  is  horizontal subspaces   $H \subset  T_pP$  such  that   $TpP = T^v_pP \oplus H$.\\
        The   differential  $d \theta$ of the   soldering  1-form   defines  a function  $\tau:  J^1 \to  \mathrm{Tor}(V)$
        with   values   in  the  space   $\mathrm{Tor}(V):= V \ot \Lambda^2(V^*)$ of $V$-valued  2-forms.
     It is  called  the {\bf  torsion  function} and it  associates   with  a horizontal  space $H= H_p $  the  2-form  $\tau_H$ defined  by
     $$    \tau_H(u,v) = d \theta (u^H, v^H) \in  V,\,\,  u,v \in V$$
  where  $u^H, v^H$ are the horizontal  lifts to $H \subset T_pP$  of  tangent vectors $\hat{p}u, \hat{p}v \in  T_{\pi(p)}Q$.\\
  \subsubsection{  C-geodesics   of    a   $G$-structure of  type $k=0$}
   Assume    that   the linear  Lie  algebra  $\gg = \mathrm{Lie}(G) \subset \ggl(V)$ has type $k=0$, i.e.  has trivial  first prolongation
   $$   \gg^{(1)} := \gg \otimes V^* \cap V \ot S^2(V^*) =0.$$
    Then   the Spencer  differential
    $$ \p : \gg\ot V^* \to  \mathrm{Tor}(V)$$
     $$      \p (A \ot \xi)(u,v) = \xi(u)Av - \xi(v)Au \in  V         $$
    is  an  embedding.   Assume  that  there   is  a  $G$-invariant  complementary subspace  $W$  to  the  image
    $\p(\gg \ot V^*)$   in $\mathrm{Tor}(V)$.   Then the preimage    $ D: = \tau^{-1}(W)$,  where  $\tau$ is  the  torsion  function, see  4.3.1,
     is   a $G$-invariant  distribution of horizontal  subspaces. More precisely, for  any $p \in P$  there  is a unique horizontal  subspace  $D=D_p$  such  that
    $\tau_H \in W$ and  the  field  $ P \ni p \to  H_p$ is  $G$-invariant.
    Such  distribution defines  a linear  connection in the frame  bundle  $ \pi :P \to Q$ with  the  connection  form
     $\omega  : TP =  T^vP \oplus  D \to \gg$  which  has  kernel $\ker \o =  D$ and   coincides  with the vertical parallelism  $T^vP \to \gg$ on    $T^vP$.  Then   the   sum
    $$\kappa =  \theta + \o : TP \to V \oplus \gg$$
    is  a Cartan  connection.
    Like  in  the Riemannian case, {\bf  C-geodesics} are  defined  as  projection  to  $Q$ of constant  horizontal vector  fields $X \in \kappa^{-1}(V)$, and they coincides with    geodesics of    the linear  connection  $\o$, see \cite{K-N}.


   \bex \label{Exercise}    Let $(Q,g)$ be  a Riemannian manifold  and   $ \pi : P\to Q$ the $ O(V)$-bundle of orthonormal  frames, i.e. $ O(V)$ -structure.\\
    One  may easily  check  that   the first prolongation
        of  the  orthogonal Lie  algebra $\gso(V)$ is  trivial  and  that   the  map $\p : \gso(V) \ot V^* \to Tor(V)$ is  an isomorphism. Taking $W=0$,   we get the distribution  $D = \tau^{-1}(0)$, defined  by  the  condition  $\tau|D =0$. The  associated  connection  is  the Levi-Civita  connection.\\
 More generally, let  $G \subset O(V) $ be  a   closed  subgroup of  the  orthogonal group  $O(V)$ and $\pi : P \to Q$ a  $G$-structure.  Let  $\gg^{\perp}$ be the orthogonal  complement   to   the   subalgebra   $\gg = \mathrm{Lie}(G)$
  in   $\gso(V)$ with  respect  to  the Killing  form.   Then $W = \p(\gg^{\perp} \ot V^*)$ is a $G$-invariant  complement  to   $\p(\gg\ot V^*)$  in  $\mathrm{Tor}(V)$. The   corresponding  distribution  $D = \tau^{-1}(W)$  defines a linear  connection $\o^{can}$ in  $\pi$   with  torsion in $W$.
   In  the  classical language,  this means  that  at  any point $q \in Q$,  the   torsion  tensor $T_q$ of  the linear  connection   $\o$ at  a point $q \in Q$, calculated   with  respect  to   a frame $\hat{p}$,  takes  values in  the   subspace  $W = \gg^{\perp} \op V^*$.\\
\eex
 We   call  the  connection  $\o^{can} $ the {\bf  canonical  connection  of  the $G$-structure with  $G \subset O(V)$.}
  We use  this  example  for the  definition  of  the \SR geodesics in the  sense of  T. Morimoto.\\

\subsubsection{  C-geodesics   for  a   $G$-structure  of   finite  type  $k>0$}

   Assume  that  $G \subset GL(V)$  group $G$  has  the  finite  type $k$,  that  is   its  Lie  algebra  $\gg$  has  non  trivial $k$-th  prolongation  $\gg^{(k)}$  and  $\gg^{(k+1)} =0$.  Then   the  full prolongation
   $$ \gg^{(\infty)} =\sum_{j=-1}^{\infty} \gg^{(j)}= V \oplus \gg \oplus \gg^{(1)} \oplus \cdots \oplus \gg^{(k)},\,   V =\gg^{(-1)},\,  \gg = \gg^{(0)} $$
    is  a  finite  dimensional   $\bZ$-graded Lie  algebra, see \cite{Stern}.
    The bundle   $\pi$ can be prolonged (\cite{Stern})  to  a bundle $ \pi^{(k)}: P^{(k)} \to  Q$   with  absolute parallelism
 $$\kappa : TP^{(k)} \to \gg^{(\infty)} = V \oplus \gg \oplus \gg^{(1)} \oplus \cdots \oplus \gg^{(k)}.$$
 {\bf C-geodesics for  the  $G$-structure  $\pi : P \to Q$ of  finite  type   $k$ }  are  defined   as  the projection  of  orbits   of  constant  vector  fields   $ X \in \kappa^{-1}(V) \subset \cx(P^{(k)})$ to $Q$.\\
\br    In general, $ \pi^{(k)}: P^{(k)} \to  Q$   is not  a principal bundle  and  $\kappa$ is not  a Cartan  connection.
\er
Of particular interest is the case  of $G$-structures,  when   $G \subset GL(V)$ is  an irreducible linear Lie  group of type  $k=1$. Then  the  full  prolongation
$$\gg^{\infty} =  \gg^{-1} \oplus \gg^{(0)} \oplus \gg^{(1)}  =  V \oplus  \gg \oplus V^*$$
is a simple 3-graded Lie  algebra.  List of  all  such 3-graded Lie  algebras  is known  and  it is very  short.
  In  this case, the  prolongation  $ \pi^{(1)}: P^{(1)} \to Q $ is  a principal  bundle   with  the  structure group
  $G^{\geq 0} =  G\cdot G^{(1)}$,  associated  to  the Lie  algebra  $\gg^{\geq 0} = \gg + \gg^{(1)}$.
 Moreover,  there  exists  a  canonical  choice  of the  absolute  parallelism $\kappa : T P^{(1)} \to \gg + \gg^{(1)}$
 ( the  {\bf normal Cartan  connection} )   which is    $G^{\geq 0}$-equivariant  and has  the   coclosed curvature, see \cite{K}, \cite{C-S}, \cite{A-D}, \cite{A-M-S}.\\
    The   C-geodesics for   such   geometries  form  an interesting  class  of  distinguished   curves in  $Q$, studied, for example, in \cite{C-S}, \cite{C-S-Z}, \cite{D-Z}, \cite{H}.\\
     For   the  conformal  structure,  which can   be  considered  as
    $\bR^+ \cdot SO_n $-structure,    generalized   geodesics  are  conformal  circles.\\

\section{  Regular  sub-Riemannian manifolds   and   Morimoto  geodesics}

     Now   we  consider   less  familiar    approach  for the   definition  and the study of \SR  geodesics, based on the   theory of Cartan
      connections  and Tanaka  theory  of  non-holonomic  $G$-structures (or Tanaka  structures).

\subsection{Regular  distributions  and  regular  sub-Riemannian  structures}
  \subsubsection{Regular distributions}
   Let   $D \subset TM$ be  a  bracket generating  distribution on  $M$  and    $  \cd^{-1}:= \G \cd$   the $C^{\infty}(M)$-module  of  sections. It  generates   a  negative filtration
   $$ \cd^{-1}\subset \cd^{-2}\subset \cdots \subset \cd^{-k} =  \cx(Q)  $$
    of   the  Lie  algebra of  vector  fields,  inductively defined by
       $$\cd^{-i-1} := \cd^{-i} + [\cd^{-i}, \cd^{-i}], \,\,  i =1, 2, \cdots.$$
   The   restriction  $D^{-i}_q = \cd^{-i}|_q$ of vector  fields  to a point $q\in Q$  defines  a  flag  of    subspaces
   \begin{equation}\label{flagofsubspaces}
   D^{-1}_q \subset D^{-2}_q \subset \cdots \subset D^{-k}_q = T_qQ
\end{equation}
   of  the tangent  space.
    The  associated  graded  space
    $$   T_q^{gr}Q =  \gm_q = \gm_q^{-1} \oplus \gm_q^{-2} \oplus \cdots \oplus \gm_q^{-k} := D_q^{-1} \oplus D_q^{-2}/D_q^{-1} \oplus \cdots +\oplus D_q^{-k}/D_q^{-(k-1)} $$
has  the  structure of   negatively graded  metric Lie  algebra, induced  by  the Lie bracket of vector  fields.
  The  graded Lie  algebra   $\gm_q$   is called  the  {\bf  symbol algebra of  the   distribution $D$ at  a point $q$}
  or  {\bf graded tangent  space   at   $q$.}\\
   The distribution  $ D $ is  called     a   {\bf  regular  distribution  of  type   $\gm$ and  depth  $k$},
    if all  symbol  algebras   $\gm_q, \,  q\in Q$ are isomorphic  to  a  fixed negatively graded Lie  algebra
    $\gm = \gm^{-1} \oplus \cdots \oplus \gm^{-k}.$
      Then  (\ref{flagofsubspaces})   defines   the  {\bf  derived  flag }   of vector  bundles
       $$D^{-1}=D \subset D^{-2} \subset \cdots \subset D^{-k} = TQ.$$
         Note  that   $\gm$ is   a {\bf fundamentally graded Lie  algebra}, i.e. it is   generated  by   $\gm^{-1}$.

\subsubsection{Regular  sub-Riemannian   manifolds } 

  Let  $(Q, D, g)$ be  a \SR manifold, where  $D$ is a regular  distribution of  type  $\gm$.

    Then    the   graded  tangent   space   $T_q^{gr}Q = \gm_q$  has  the  structure   of a {\bf   negatively  graded metric Lie   algebra},     i.e.   a  graded  Lie  algebra $\gm_q = \sum_{i=-1}^{-k} \gm_q^i$  with  an  Euclidean metric $g_q^{\gm}$   such  that  the  graded  spaces $\gm_q^i$  are mutually orthogonal.\\

The  metric   $g_q^{\gm}$  is  a  natural extension  of  the  \SR metric  $g_q^D$ in $  D_q$,  which is  described in   the  following elementary  lemma.
\bl
 Let  $\gm = \gm^{-1} + \cdots  + \gm^{-k}$ be   a   negatively graded fundamental  Lie  algebra. Then   an Euclidean metric  $g$  on  $ \gm^{-1}$   has  a natural   extension   to   an Euclidean metric $g^{\gm}$ in  $\gm$.
  \el

 A \SR manifold  $(Q,D,g^D)$  with  a  regular distribution $D$ of  type  $\gm$ is  called  a {\bf regular \SR manifold  of   the  metric type  $(\gm,  g^{\gm})$ } if all    metric  Lie  algebras  $(\gm_q, g_q^{\gm}) $  are  isomorphic  to  the metric  graded Lie  algebra  $(\gm,  g^{\gm})$.\\

\subsubsection{Regular  sub-Riemannian  structure as Tanaka  structure}

 Let  $D\subset  TQ$ be  a  regular  rank-m distribution  of  type $\gm$,    and
 $\mathrm{Aut}(\gm)$ the  group of   graded   preserving   automorphisms  of  $\gm$.\\

 An  {\bf admissible  frame}  of  $D$  is  an  isomorphism
 $$f : \gm \to T^{gr}_qQ = \gm_q $$
  of  graded Lie  algebras.  The automorphism  group  $\mathrm{Aut(\gm)}$ acts  freely  and properly  on the   manifold
   $\mathrm{Fr}(D)$ of  admissible  frames on  $D$   with  the  orbit  space  $\mathrm{Fr}(D)/ \mathrm{Aut}(\gm) =Q$.  Hence,  $\mathrm{Fr}(D) \to  Q$ is   a principal  bundle (  called {\bf the bundle of   admissible  frames on  $D$}).\\
   Let   $G^0 \subset \mathrm{Aut}(\gm)$   be  a  Lie  subgroup.   A  {\bf Tanaka   $G^0$-structure} ( or   {\bf a relative $G^0$-structure})  is  a   $G^0$-principal  subbundle    $\pi : P \to  Q = P/G^0$  of the bundle of  admissible  frames  on  $D$.\\

The      classical  identification    of    Riemannian manifolds    with   $O_n$-structures is  extended  to the   \SR case:\\
 \bp A    regular sub-Riemannian   manifold $(Q,D, g^D)$ of    type   $(\gm, g^{\gm})$ is   identified  with  a Tanaka  $G^0$-structure  with  the  structure  group  $G^0 = \mathrm{Aut}(\gm, g^{\gm}) \subset O(\gm)$. Conversely, any $G^0=  \mathrm{Aut}(\gm, g^{\gm})$-Tanaka structure  defines a regular sub-Riemannian   manifold $(Q,D, g^D)$ of    type   $(\gm, g^{\gm})$.
\ep
\pf    The  Tanaka $G^0$-structure,  associated  to $(Q,D,g^D)$, consists   of  all admissible frames $f : \gm \to \gm_q$,  which    are  isomorphisms  of    the metric  Lie  algebras.  \\
Conversely, let $(\pi : P \to Q, \, D)$
 be   a Tanaka $G^0$-structure with  $G^0= \mathrm{Aut}(\gm, g^{\gm})$. Then the  associated  sub-Riemannian  metric on  $D$ is   defined    by  the  condition   that   for   any admissible  frame
          $f \in  P$, its      restriction
          $$f_{\gm^{-1}} : \gm^{-1} \to  D_q = \gm^{-1}_q$$
            is  an  isometry. \qed\\

\subsection{Morimoto  geodesics   of sub-Riemannian manifolds}

\subsubsection{\label{Tanaka-Morimoto} Tanaka   prolongation  of    non-holonomic  $G$-structures }  

N. Tanaka  generalised  the  theory of  $G$-structures to  Tanaka  structures.   In  particular,  he   defined    the   full prolongation   of  a  non-positively graded Lie  algebra $\gg =  \sum_{i=-k}^{0}\gg^i = \gg^{-k} +  + \gg^{-1} + \gg^0$   as   a  maximal   $\bZ$-graded   Lie algebra
 of  the  form
   $$ \gg^{(\infty)}  = \sum_{i=-k}^{\infty}\gg^i =  \gg \oplus \gg^{(1) }\oplus \gg^{(2)} \oplus \cdots$$
  such  that  for  any   $X \in  \gg^{(i)}, \,  i>0$  the  condition  $[X, \gg^{-1}]=0$ implies   $X=0$.\\
   A non-positively  graded Lie algebra $\gg$ is  called a {\bf  Lie  algebra  of  finite  type  $\ell$} , if $\ell < \infty$  is   the maximal number   such  that   $\gg^{(\ell)} \neq 0$.
 \bt (Tanaka)(see \cite{C-S},\cite{T}, \cite{Z},\cite{A-D1}).   Let   $\pi:  P \to Q$   be  a Tanaka  $G^0$-structure   on $(Q,D)$   where  $D$   is  a  regular   distribution  of  type  $\gm = \gm^{-1} + \cdots + \gm^{-k}$    and $G^0 \subset \mathrm{Aut(\gm)} $ a connected  closed  subgroup of  the   automorphism group   with Lie  algebra
 $\gg^0 \subset  \mathfrak{der}(\gm)$.  Assume  that
   the non-positively graded  Lie  algebra  $\tilde{\gm}:= \gm + \gg^0$ has finite  type   $\ell$.
      Then there is  a  canonical bundle
    $  P^{\infty}  \to  Q$,  constructed  by  successive prolongations, with  an  absolute  parallelism
     $\kappa :  T P^{\infty} \to \tilde{\gm}^{\infty}$.
If   the  first prolongation    $\gg^{(1)} =0$,  then  the  absolute parallelism   $\kappa :  TP  \to  \tilde{\gm}= \gm + \gg^0$  is    a Cartan  connection.
\et
\subsubsection{ Morimoto definition of sub-Riemannian geodesics}
The  Morimoto definition  of sub-Riemannian geodesics  is based  on the  following
    important  theorem ,  which  he proved in  the framework of  his remarkable  theory of filtered  manifolds \cite{M1}.


\bt (T. Morimoto \cite{M}) i) Let  $(\gm = \sum_{-k}^{-1} \gm^i,  g^{\gm})$ be  a  fundamental  negatively graded  metric  Lie  algebra   and    $G^0 = \mathrm{Aut}(\gm, g^{\gm})$    the linear Lie  group  of  orthogonal  automorphisms    with  the  Lie algebra  $\gg^0 = \mathfrak{der}(\gm, g^{\gm})$.
Then  the   full prolongation  of  the non-positively graded Lie  algebra
$\tilde{\gm} = \gm + \gg^0$   coincides  with  $\tilde{\gm}$.\\
ii)   Let  $(Q, D, g^{D})$   be   a   regular  sub-Riemannian  manifold  of the  metric  type  $(\gm, g^{\gm})$.
Then  the  associated  Tanaka  structure    $\pi: P \to Q = P/G^0$    admits  a canonically  defined 1-form
\be \label{Morimotoform} \kappa :  TP \to   \tilde{\gm} = \gm + \gg^0
\ee
   such  that  $(\pi, \kappa)$  is  a normal  Cartan  connection of  type   $L/G^0$, where $L$ is the  simply  connected Lie  group   associated  with
   the Lie  algebra  $\tilde{\gm}$.
 Moreover, the  horizontal part    $ \theta = \mathrm{pr}_{\gm} \circ \kappa$  of $\kappa$ is  a  soldering form  and  the  vertical part
 $\o = \mathrm{pr}_{\gg} \circ \kappa$ is  a principal  connection.
 \et

   The {\bf  Morimoto  sub-Riemannian geodesics} (  shortly, {\bf  M-geodesics}) of  a regular \SR manifold $((Q,D,g^D)$ is  defined   as  projection to  $Q$ of the integral  curves  of  a constant vector  fields $X \in   \kappa^{-1}(\gm^{-1})$,  where   $\kappa$ is   the associated
    Cartan connection  (\ref{Morimotoform}).

 \subsection{ Admissible  rigging,  associated Cartan  connections  and   Cartan-Morimoto  geodesics }

 Here  we develop    an  elementary   approach   for   constructing   Cartan  connections (in particular, normal Cartan connection ) for
  regular  \SR  manifold, considered  as Tanaka structurs.   It is   working  also  for   other Tanaka  structures  with  trivial  first prolongation. It is based  on  the notion  of  admissible rigging,  see \cite{A-M-S},  and    results  from  \cite{A-D}.\\
 Let   $(Q, D, g^D)$ be  a regular \SR manifold   of    metric type  $(\gm, g^{\gm}) $,  and
  $ D^{-1} = D \subset D^{-2} \subset \cdots \subset D^{-k}=TQ  $
  the  derived   flag   of  distributions. \\
 A complementary to  $D$  distribution  $V$  with  a direct  sum  decomposition  $V = V^{-2} \oplus  \cdots \oplus V^{-k}$
 as called   an  {\bf admissible  rigging} if  $V^{-j}$  as  a complementary  to $ D^{-j+1} $ subdistribution in $D^{-j}$.
 In other  words,
 $$ TQ =  D \oplus V^{-2} \oplus \cdots \oplus V^{-k},  \, D^{-j } = D^{-j+1} \oplus V^{-j}, \,  j=2,\cdots,k.
    $$
   Since   $\gm^{-j} :=D^{-j}/ D^{-j+1}= (D^{-j+1}\oplus V^{-j})/D^{-j} \simeq  V^{-j} $,  the  admissible  rigging $V$
 defines  an  isomorphism  $\psi_V : T^{gr}Q  \to  TQ   $ of  vector    bundles.
   It   induces  an  isomorphism  $\hat{\psi}_V : P \to P^V$ of  the Tanaka  structure   $ \pi: P \to Q = P/G^0$, associated  to  the \SR manifold,    onto  a $G^0$-structure,  which   we denote  by
     $\pi^V : P^V \to  P^V/G^0$. Identifying   these  principal bundles,  we will  consider   the  soldering  form $\theta^V$ of  the $G^0$-structure  $\pi^V$   as   a  soldering   form  on  $P$.
It  turns the Tanaka  principal bundle
  $\pi:  P  \to Q $ into a  $G$-structure,  with   the  structure  group  $G^0 = \mathrm{Aut}(\gm, g^{\gm}) \subset O(\gm)$. We  denote  by $g$   the associated  Riemannian metric in  $Q$  and  by     $\o^V$  the  canonical  connection  form.
   The  principal   connection  $\o^V$   defines a Riemannian  connection $\nabla^{\o^{V}}$  with  torsion. Denote  by
    $\tau_{\theta^V} : P \to  W := \p(\gg^{\perp} \ot \gm^*)   \subset Tor(\gm)  $  the  associated  torsion  function.
   Then  the 1-form
$$\kappa^V :=  \theta^V + \o^V:   TQ  \to  \gm +\gg $$
 defines   a  structure  of Cartan  connection  in  the  principal bundle  $\pi : P \to  Q$.\\
This  Cartan  connection  induces     the  Tanaka  structure  $\pi: P \to Q$  via   the
isomorphism    $\psi :  TQ \to T^{gr}Q$.  Moreover, any   Cartan  connection, which induces  the  Tanaka  structure  $\pi$, is  associated  with  some  admissible rigging,  see  \cite{A-D}, proposition 2.\\
   Like  in   the case  of  normal  connections,  we define \SR geodesics   as   the   projection  to  $Q$ of  the integral  curves  of  constant vector  fields  from  $\kappa^{-1}(\gm^{-1})$. We call  such geodesics  {\bf CM-geodesics   associated  to  an admissible  rigging}.
   Since the   connection  $\n^{\o^{V}}$ preserves  the   \SR metric  $(D, g^D)$,  CM-geodesics  are  $D$-horizontal geodesics
   of  the  connection $\n^{\o^V}$  or, in other  words, the  geodesics of  the partial  connection  on $D$, which is  the  restriction of
      $\n^{\o^V}$  to $D$.  Note  that    this partial  connection  coincides  with the Schouten partial connection  $\n^S$ associated to  the  rigging $V$  if  and only if  the  torsion  function $\tau$ of the soldering  form $\theta^V$
      satisfies  the  condition   $\tau_p (D,D) \subset V, \, \forall p \in  P$. We  get

 \bp
     Let  $V$ be  an  admissible rigging  of a regular sub-Riemannian  manifold $(Q,D,g^D)$.  Then CM-geodesics  of   the Cartan   connection  $\kappa^V$   coincide  with  the S-geodesics  of  the Schouten partial  connection  $\n^S$,  defined by the rigging $V$,                                                  
      if and   only if  the  torsion  function $\tau $
        of the soldering  form $\theta^V$
      satisfies  the  condition   $\tau_p (D,D) \subset V, \, \forall p \in  P$.
 \ep

\subsubsection{  Admissible  riggings  and    the normal Cartan  connection }

   Here   we   apply    the  results from  \cite{A-D},  to   prove  the  existence  and  the uniqueness     of an
   admissible  rigging $V$ which  defines the  normal  Cartan  connection  $\kappa^V$ associated  to  a regular  \SR  manifold. It
     reduces  the  problem of  constructing  normal Cartan  connections  to an  appropriate   deformation  of  an  admissible  rigging.\\

   The   following  theorem  is an  elaboration   of  the Morimoto theorem.
\bt  Let  $\pi : P \to Q$ be   the Tanaka  $G^0$-structure,   associated  to a regular sub-Riemannian   manifold  $(Q,D, g^{D})$
with a metric  symbol  $(\gm, g^{\gm})$.
Then   there  is  a uniquely defined   admissible rigging $V_0$,  which  defines   a normal  Cartan  connection
$(\pi : P = P^V \to  Q, \kappa^{V_0} = \theta^{V_0} + \o^{V_0} )$.
\et

\pf  The  proof   follows  from  \cite{A-D}, theorems 1 and 2.   For the  uniqueness   of the Cartan  connection,  we  have  to  check
 that   the  first  cohomology  group for  the  cocycles  of  positive  degree vanishes:   $H^1(\gm, \gm + \gg)_1`=0 $.
 The  degree  of a cocycle  from   $(\gm^{-i})^* \ot \gm^{-j}$ is defined   as   $i-j$. The  vanishing of   this  cohomology
 is  a simple  exercise.
\qed

\bc Let $(Q,D, g^D)$ be  a regular  sub-Riemannian manifold   with metric  symbol algebra  $(\gm, g^{\gm})$  and
  $\pi: P \to Q$  the associated  Tanaka  structure  with  the  structure  group $G^0 = \mathrm{Aut}(\gm, g^{\gm} )\subset  O(\gm)$. Then \\
  i) there  exists the unique   soldering   form   $\theta^{0} : TP \to \gm$    which  together  with   the
  associated   canonical  connection  $\o^{can} : TP \to \gg$
   defines   the  normal Cartan  connection
  $\kappa^{can} = \theta^{0} + \o^{can}$.\\
ii) The sub-Riemannian  metric $g^D$ admits  a canonical  extension  to  the Riemannain metric $g$ on $Q$,  defined  by  the  condition that   the isomorphism (coframes)  $\theta^{0}_p : T_{\pi(p)}Q \to \gm$ is  an   isometry onto  the Euclidean  space $(\gm, g^{\gm})$ for $p \in P$.\\
iii) The  isometry  group $ A=   \mathrm{Iso}(Q, D,g^{D})$ is  a Lie  group of  dimension $\dim A \leq n + \dim G \leq n + \frac{m(m-1)}{2}$  where $m = \dim \gm$,  which preserves   the  Riemannian metric $g$ and   acts freely  with  closed orbits on  $P$. The  stability  subgroup $A_q$ of  a point $q \in Q$  has   the  exact  isotropy  representation
$j : A_q \to  GL(T_qQ)$  and  the isotropy  group  $j(A_q)$  is identified  with  a  subgroup of  the   group
  $G^0 = \mathrm{Aut}(\gm, g^{\gm})$ and  has  dimension  $\dim A_q \leq \frac{m(m-1)}{2}$.
\ec

\pf i)  and ii) directly  follow  from the Theorem. iii) The  isometry  group  $A$ preserves   the   absolute  parallelism
 $\kappa^{can}$ on  $P$ and  the Riemannian metric  $g$. By the Kobayashi theorem about  the  automorphism  group  of  an  absolute parallelism \cite{Stern}, it  acts  freely    with  closed orbits  on $P$. Hence $\dim A \leq \dim P = n +  \dim G \leq \frac{m (m-1)}{2}$,   since   $G$ is identified  with   a  subgroup of  the orthogonal group $O(\gm)$.  The  stability  subgroup $A_q$ has exact isotropy  representation
 $j : A_q \to  GL(T_qQ) $  and  since $j(A_q)$  preserves  the metric of $T_qQ$ and  the  derived  flag
 $ D_q \subset D^{-2}_q \subset \cdots \subset D^{-k}_q$, it  acts  freely in  $T_qQ \simeq T^{gr}Q = \gm_q$  and  preserves  the  Lie  algebra structure  in    $\gm_q$. This   implies  that   $j(A_q)$ is identified  with  a  subgroup  of
$G$.
\qed


\section{  Geodesics  of  Chaplygin  transversally \\
 homogeneous  sub-Riemannian manifolds}  
 \subsection{Principal  connection  and  its  curvature}
 Let $\pi:  Q \to  M = Q/G $ be  a  $G$-principal  bundle  with     a right  action  $R_g q = qg$    of a  Lie  group  $G$.\\
 For  $a \in  \gg = \mathrm{Lie(G)}$, we  denote  by
     $a^* : q \mapsto  qa :=  \frac{d}{dt} q \exp(ta)|_{t=0} \in T_qQ$
    the fundamental  vector  field  (  the  velocity vector  field of $R_{\exp ta}$).\\
 Recall  that  the principal  connection is a  $G$-equivariant   $\gg$-valued   1-form  $\vpi: TQ \to \gg$,  which
   is  an  extension  of  the   vertical parallelism,  defined  by  $ T_qQ \ni a^*_q \mapsto a \in \gg$.
  The  equivariancy means  that
  $               R_g^* \vpi = \Ad_g^{-1}\circ \vpi,\,  g \in G.$

The  connection  $\vpi$   is completely    determined  by   the  horizontal $G$-invariant distribution   $D = \ker \varpi$.\\
 Since   $T_qQ =  D_q \oplus T^v_qQ$,  any   vector  $X \in T_qQ$  is  decomposed as   $X = X^h \oplus X^v$ into the horizontal  and  the vertical parts.
   Moreover,   any vector   filed   $X \in \cx(M)$   has  the  canonical ($G$-invariant)  horizontal lift   $X^D \in  \G D$
       due  to   the     isomorphism   $\pi_* : D_q \to  T_{\pi(q)}M$. Recall  the  following  formulas,   \cite{K-N},
  \be
  \label{K-N  formula}
   [a^*, b^*]  =  [a,b]^*, [a^*, X^D] =0, \, [X^D,Y^D]=  [X,Y]^D +  [X^D,Y^D]^v .
  \ee
where  $a,b \in \gg, \,  X,Y \in \cx(M)$.\\
Note   that  $[X^D,Y^D]^v$  is  a  $G$-invariant   vertical  vector  field  and  its  restriction  to a  fiber
$ \pi^{-1}(x) $ depends only  on  vectors  $X_x,Y_x$. So   we have  a linear  skew-symmetric   map
$$  \ca : T_xM \times T_xM \to \cx(\pi^{-1}(x))^G, \,\, (X,Y)\mapsto \ca(X,Y) := \frac12 [X^D,Y^D]^v $$
 into  the  space  of  $R_G$-invariant vector  fields   along  the  fiber  $\pi^{-1}(x)$. \\
  For  $X \in \cx(M)$,   we    defines a $C^{\infty}(Q)$  linear map
    $$   \ca_X : \G (D) \to \G T^vQ, \, Y^D \mapsto  [X^D,Y^D]^v =  \ca(X,\pi_* Y^D)  $$
     which maps  $G$-invariant  horizontal vector  fields $Y^D$     into  $G$-invariant  vertical vector  fields from
      $\G T^v(Q)$.
     The  dual   map
    $$\ca^*_X : \G(T^v(Q)) \to  \G  D $$   sends   $G$ -invariant vertical vector  fields  into  $G$-invariant  horizontal vector  fields.\\
 The  curvature of   $\vpi$ is  the $\gg$-valued   horizontal 2-form  $F =  d \vpi + \frac12  \vpi \wedge \vpi$.  Denote  by
 $e_{\a}^*$   the   basis   of  fundamental vector  fields,  associated  to  a basis  $e_{\a}$  of  $\gg$.  Then the  curvature  form  is  related  with the  tensor  $\ca(X,Y)$  as  follows:

  \be \label{curvature}
  \begin{array}{lll}
    F_q(X^D,Y^D) &=  F_q(X^D,Y^D)^{\a}e_{\a}  & = -\vpi_q([X^D,Y^D]) + (\vpi_q \wedge \vpi_q )(X^D,Y^D)\\
               & = -2 \vpi_q (\ca_q(X, Y))&= -2\vpi_q( \ca_q^{\a}(X,Y)e^*_{\a})\\
               & =    (\ca_q(X, Y)^{\a}  e_{\a}.
   \end{array}
  \ee

\subsubsection{Formulas  in   coordinates   }
 To  write  formulas  in  coordinates,  we    fix   a  section  $s$  of $\pi$.  To  simplify notations,  we  assume  that
 the  section  $s :  M \to  Q$ is  global.  It  defines    a  trivialization
$$  M \times G   = Q , \, \,\,     (x,g) =  s(x)g
   $$
of  the  principal  bundle,  where the  group   $G$ acts on  $M \times  G$  as
$$   R_g (x,g_1) = (x, g_1 g).
$$

   A fundamental   vector  field  $a^*,  \,  a \in \gg$   is identified    with   the  left invariant  vector   fields
     $$a^L : (x,g) \mapsto  (x, ga) := (L_g)_* a.$$
  We  denote    by   $a^R : (x,g)\mapsto (x, ag)$
    the  right invariant vector  field,  generated  by  $a \in \gg$.\\
We  identify   the   tangent   bundle  $TG $ with   the   trivial bundle  $\gg \times G $   (  which  we  denote  simply  as $\gg$), using left translations
$$        \gg \times G \ni (\dot g, g)   = g^{-1}\dot g := (L_g)_*^{-1}\dot g  \in  T_gG  $$
 and  the  tangent  bundle  $T(M \times  G) = TM \times  TG$  \, with  the  bundle $ TM \otimes \gg$ over   $M \times G$.
  The  tangent  vector  $ \dot x + ga :=   \dot x + (L_g)^*a \in  T_{(x,g)}(M \times G)$  will be   denoted  by   $(\dot x,  a)$, where  $\dot  x \in  T_xM, \,  a \in \gg$ .\\
 The pull back $\o^s := s^* \vpi$  of  the   connection   form $\vpi$  to  $M$ is  a $\gg$-valued  1-form
  $\o^s = A^{\a}_i (x)dx^i \otimes e_{\a}  \in \Omega^1(M, \gg)$  on  $M$.  The horizontal  distribution  along  $M \times \{e\}$ is  given  by
  $$D_{(x,e)}   = \{ \dot x - \o^s(\dot{ x}) \} =   \{\p_i - A^{\a}_i (x) e_{\a}\}. $$
  Since the  distribution $D$  is  $R_G$-invariant,  we  get
  $$    D_{(x,g)} =  (R_g)_* D_{(x,e)} = \{ \dot x  - A^{\a}_i(x)e_{\a}^R\} =
  \{\dot x  - A^{\a}_i(x) (\Ad_g)_{\a}^{\beta}e_{\beta}^L\}        $$

    Note  that $\vpi(e_{\a}^*) \equiv \vpi(e_{\a}^L) = e_{\a}  $.  This  shows   that $\vpi|_{TG}$   coincides  with    the  left invariant   Maurer-Cartan  form $\mu (\dot g) =  g^{-1} \dot g  $. Hence,  the  connection   form may be  written  as  $\vpi = \mu  + {A}$,  where
     ${A} = A^{\a}_i (x,g) dx^i \otimes e_{\a}$ is  a  1-form, which vanishes on  $TG$. Solving the  equations
 $$  0=  \vpi (\p_i^D)= \vpi(\p_i -A^{\a}(x)e_{\a})=- \mu(A_i^{\a}(x)e_{\a} ) +{A}^{\a}_i(x,g) e_{\a}=$$
 $$
  -A^{\a}_i(x) (Ad_{g}^{-1})^{\beta}_{\a} e_{\beta})  + {A}^{\a}_i(x,g) e_{\a},$$

 we   find  that  the  connection  form $\vpi$ on   $Q = M \times G$ is  given  by

  \be \label{connectionform}
  \vpi = \mu + A  ,\,\,   A  =  A^{\a}_i (x,g)dx^i \ot e_{\a} =
    A^{\a}_i(x)  dx^i \ot (\Ad_g^{-1})^{\beta}_{\a} e_{\beta}.
  \ee

  \subsection{ Chaplygin  metric  and  its  standard  extension   }

   Let $  \pi : Q \to M = Q/G$ be  a principal bundle  with  a  connection $\vpi$   and  $D = \ker \vpi$ the  horizontal  distribution.

   A Riemannian  metric  $g^M$  on   the  base manifold  $M$ defines  a canonical  invariant  \SR  metric   $g^D$  on  $D$
       such  that   the  projection $\pi_* :  D_q \to  T_{\pi(q)}M$ is  an isometry.

 The   \SR   metric   $(D, g^D)$ is  called  a Chaplygin metric and    the  \SR  manifold   $(Q,D,g^D)$  is called  a  {\bf  Chaplygin  system}
   or a {\bf transversally    homogeneous  \SR  manifold.}

 Let   $(Q,D,g^D)$    be  a Chaplygin  system  associated  to  a principal bundle $(\pi :  Q \to  M,  \vpi )$  over  a Riemannian manifolds   $(M,g^M)$  as  above  and  $g^{\gg}$
   an Euclidean metric  on the  Lie  algebra  $\gg$.
It  defines  a degenerate  metric
         $$    g^F(X,Y) = g^{\gg}(\vpi(X), \vpi(Y))$$
      on  $Q$     with kernel  $D$,  whose  restriction   to a  fiber
        $\cf(x) =  \pi^{-1}(x)$  is  a  Riemannian metric. We   will  consider  also  $g^D$ as  a degenerate metric  on  $Q$  with  $\ker  g^D = T^vQ$.

       Then
       \be \label{metric in LA}   g^Q =  g^F  \oplus  g^D
       \ee
     is  a  Riemannian metric  in  $Q$. It  is  called   the {\bf standard  extension  of   the Chaplygin  metric  $g^D$. } \\
      Note  that    the   metric   $g^Q$ is     $G$-invariant  if and  only  if the    degenerate  metric $g^F$ is invariant or,  equivalently,   the metric  $g^{\gg}$  is  $\mathrm{Ad}_G$-invariant. In  this case  the   metric $g^Q$ is called  the {\bf  bi-invariant  extension } of  the   \SR metric  $g^D$.\\
 Denote by   $\n^M$  (resp., $\n^Q$) the  Levi-Civita  connection  of $g^M$
 (resp.,  $g^Q$),  and by $\n^F$    the  Levi-Civita  connection  of  the  induced   metric  $g^F$ on a  fiber,  which is  a totally geodesic  submanifold of  $(Q,g^Q)$.\\
The  Koszul formula   implies  the  following  O'Neill  formulas  for  the  covariant  derivative  of   fundamental  vector  field $b^*$ and  the  horizontal lift $X^D$  of  a vector  field $X \in \cx(M)$, see   \cite{Besse}.\\

 $$
 \begin{array}{lll}
    i)&  \n^Q_{a^*}b^*   =&  \n^F_{a^*}b^*, \, \,\,   a,b \in \gg\\
   ii)&  \n^Q_{a^*}X^D = & \n^Q_{X^D}a^*  = (\n^Q_{a^*}X^D)^h = \ca^*_{X}a^*, \\
  iii)& \n^Q_{X^D}Y^D =& (\n^M_X Y)^D   + \ca(X,Y),\, \,  X,Y  \in  \cx(M).
 \end{array}
 $$
 The connection $\n^F$ is described  in  terms of Lie  brackets  as  follows
   $$ 2g(\n^F_{a^*}b^*, c^*) = g^{\gg}([a,b],c)- g^{\gg}(b, [a,c])- g^{\gg}(a,[b,c]),\,  a,b,c \in \gg. $$

  \subsubsection{   S-geodesics of  a Chaplygin metric}

The O'Neill formulas   imply  the  following  relations  between  S-geodesics  of  the Chaplygin   \SR metric  and  geodesics  of  the Riemannian metrics  $g^Q, g^M$, see
also \cite{Besse}.
\bt \label{S-geodesics}  i) The principal bundle  $\pi  : Q \to  M$   with a  standard   metric $g^Q$  associated  to a Riemannian metric $g^M$  is  a   Riemannian  submersion  with  totally geodesic  fibers.\\
ii)A Riemannian  geodesic  $\g(s)$ of   $(Q, g^Q)$  which   is  horizontal  at one point   is   horizontal  and   it projects     onto the  geodesic   $\pi \g(t)$  of the  base manifold  $(M, g^M)$.\\
iii) S-geodesics  of   the   Chaplygin metric    are precisely  horizontal geodesics  of  $(Q,g^Q)$   and  they  are  horizontal lifts of  geodesics  of  the base manifold.

\et

\pf Recall  that S-geodesics  are  geodesics of  the Schouten  connection
$\n^S_X Y = \mathrm{pr}_{D} \n^Q_X Y, \, X,Y \in  \G  D$,  that is     horizontal
curves $\g(s)$ which  satisfy   the  equation
$$            0=    \n^S_{\dot \g(s)}\dot {\g}(s)= \mathrm{pr}_{D} \n^Q_{\dot \g(s)}\dot {\g}(s)= (\n^M_{\dot x} \dot x)^D .$$
The  result  follows   from   the O'Neill  formula iii) and  skew-symmetry  of  the  tensor  $\ca(X,Y)$.\qed \\

\subsubsection{H-geodesics  of  a Chaplygin metric}
 Now   we  consider  H-geodesics   of the Chaplygin metric $g^D$ on the  principal bundle  $ (\pi: Q \to  M, \vpi)$ over  a Riemannian manifold  $(M, g^M)$   and study their relation  with  geodesics of  the
 standard extension  $g^Q$ and  the  metric $g^M$ of  the  base manifold.\\
As   above,   we  fix  a  trivialization $Q = M \times G$  of the  principal bundle  $\pi$,  defined  by a  section  $s$. Then  any  fiber  $F_x =  x \times G$ with    the  metric  $g^F$ is identified  with  the Lie  group $G$
  with  the left invariant metric (  still  denoted  by  $g^F$ )  defined  by  the  metric    $g^{\gg}$.
   Note  that   the horizontal   distribution   $D = \ker \vpi$  on  $Q =  M \times G$ and  the  \SR metric  $g^D$   are $R_G$-invariant,  but  not  $L_G$-invariant,   the fiberwise  metric  $g^F$ is  $L_G$-invariant,  but not  $R_G$-invariant  and  the  metric $g^Q$ is neither $R_G$-invariant, nor  $L_G$-invariant.\\

  The orthogonal  decomposition   $g^Q = g^D \oplus g^F$ of  the  metric  $g^Q$ defines
 the orthogonal  decomposition   $g_Q^{-1} =  g_D^{-1} \oplus  g_F^{-1}$ of  the associated  contravariant metric, which we denote  by   $g_Q^{-1}$. Here  $g_D^{-1}$ (resp., $g_F^{-1}$)is  the  contravariant   metric  on $D$ (resp.,  on $T^vQ$ ). They may be locally  written  as
 $$g_D^{-1} = \sum_{i=1}^m   X_i \ot X_i ,\,\,\,\,  g_F^{-1}=   \sum e_{\a}^* \ot e_{\a}^*$$
    where  $(X_i)$ is  a local  orthonormal  frame  in  $D$  and
 $(e_{\a})$  is an  orthonormal  basis  of  $\gg$.  We   will  consider  $g_D^{-1}$  and  $g_F^{-1}$ as functions  on  $T^*Q$ (cometrics).

       The  Hamiltonian
  $h_Q = \frac12 g_Q^{-1} \in C^{\infty}(T^*Q) $  of  the geodesic  flow   of  the  Riemannian metric  $g^Q$   is  the   sum
$h_Q =   h_D +  h_F$ of  the  Hamiltonian
 $ h_D(\xi) := \frac12 g^{-1}_D(\xi,\xi)$  of  the  \SR metric   and  the  Hamiltonian
 $ h_F(\xi) =  \frac12 g^{-1}_F (\xi,\xi)$ of  the fiberwise metric $g^F$.


 Now  we  describe  the  relations  between  \SR  H-geodesics   and   Riemannian  geodesics  of  the  metric $g^M$ on  the  base   $M$ and   left invariant  metric  $g^F$  on  the  group  $G$.   In  the  case, when  the  extension  $g^Q $ in  bi-invariant,  i.e.   the  metric    $g^{\gg}$  is   $\Ad_G$ invariant,  they had been proved  by R. Montgomery \cite{Mont}, Theorem 11.2.5 (the  main theorem).\\

\bl

   Let $g^Q = g^D  \oplus g^F  $ be a  standard  extension  of  a Chaplygin metric $g^D$.
  Then the   Hamiltonians  $h_Q, h_F, h_D$ Poisson  commute
  $$\{h_F, h_D\} = \{h_F, h_Q\}=0,$$
 and the  associated  Hamiltonian vector  fields $\vec{h}_F, \vec{h}_D, \vec{h}_F$  commute.
 \el

\pf A  fundamental   field  $a^*$
 commutes  with   the  horizontal lift  $X^D$ of  a basic  vector  field,  see  (\ref{K-N  formula}), and preserves  the
 decompositions   $TQ = T^vQ + D,\,  T^*Q = (T^vQ)^* + D^*$.   We  calculate  the  Lie  derivative  of  the   \SR metric $g^D$ as follows:
$$
 \begin{array}{ll}
    (\cl_{a^*} g^D)(X^D, Y^D)&  \\
   = a^*\cdot (g^D(X^D, Y^D))-g^D([a^*,X^D], Y^D)+ g^D(X^D, [a^*,Y^D])&\\
    = a^*\cdot g^M(X,Y)=0.&
   \end{array}
         $$
This   shows  that the fundamental field  $a^*$  preserves  the  \SR  metric  $g^D$  and  the  dual  cometric $g_D^{-1}$.
  This means  that
  $$ 0=  \cl_{a^*}g^{-1}_D  = \{ a^*, g_D^{-1}  \} \equiv \{ p_{a^*}, g^{-1}_D  \}$$
 where  $p_{a^*}$  is  the Hamiltonian of  the  fundamental  field $a^*$. Then the Leibnitz  rule  for  Poisson bracket   shows  that  $\{ a^* \ot a^*, g_D^{-1}  \} \equiv   \{ p_{a^*}^2,  g_D^{-1} \} = 0$.
 Hence the Hamiltonian   $h_F = \frac12 \sum a_{\a}^* \ot a_{\a}^* = \frac12\sum (p_{a_{\a}})^2 $ commute  with  $h_D = \frac12  g^{-1}_D$.
 \qed

  Using  the  same  arguments  as in \cite{Mont}, we   get
   \bt
    \label{H-geodesics} i)The  sub-Riemannian  geodesic flow of the sub-Rieman\-nian  met\-ric  $g^D$  is
        a  composition
  $\exp t \vec{h}_D = \exp t \vec{h}_Q\circ \exp (-t \vec{h}_F)  $
   of  the Riemannian  geodesic flows of  the  metric $g^Q$ and  the fiberwise  metric  $g^F$.\\
  ii)  Denote   by   $g_a(t) \subset   G$  the  geodesic   of  the
  group  $G$  with   the left invariant  metric $g^F$ as  above  with initial  conditions
   $g_a(0)= e, \dot{g}_a(0) =a \in \gg$
    and  by $\g_w(t)$  the  geodesic of  the standard  metric $g^Q$  with  initial  conditions $\g(0)=q \in  Q,\,  \dot{\g}(0) =w = w^v + w^h \in  T_qQ  $.
    Then  the  curve  $q(t) = \g_w(t)g_a(t)$ is a sub-Riemannian   H-geodesic   if and only if
     it has   horizontal  velocity $\dot{q}(0) = w +a^*_q = w^h\in D_q$   that  is
     $\vpi(w) = -a. $\\
   iii)  Horizontal  geodesics of $g^Q$   are  sub-Riemannian   geodesics  and   they project  to  geodesics of  $(M,g^M)$.\\
   iv) Sub-Riemannian   geodesics  are  horizontal lifts  of  the projection of geodesics $\g_w(t)$  of $g^Q$  to  $M$.
    \et

    \pf
    i) is obvious.
 ii) The  restriction  $g^F|_{F_x }$ of $g^F$ to  any  fiber  $ F_x= (x,G)$ is  identified  with the left invariant Riemannian metric (  denoted  again  by  $g^F$ ) on $G$,   defined  by the  metric  $g^{\gg}$.\\
 Note that the  projection   of  integral curves  of the  $g^Q$-geodesic  flow  to $Q$  are   geodesics  of  $g^Q$.
  The projection    of  integral  curves of $\exp t\vec{h}_F$ to  $Q$ are geodesics   of  the  metric   $g^F$, hence  also of $g^Q$,
          since the fibers  are  totally geodesics,  see O'Neill  formulas.   The  projection $q(t)=\tau\circ \xi(t)$  to $Q$ of  the composed   curves
   $$
     \xi(t) = \exp t \vec{h}_F \circ \exp t \vec{h}_Q \circ (\xi) ,\,\xi \in  T^*_qQ
$$
 are  curves   of  the  form   $ q(t)= R_{g(t)} \g(t)  = \g(t) g(t) $ where
 $\g(t) = \tau\circ \exp  t  \vec{h}_Q (\xi)$  is  a geodesic  of $g^Q$  and   $g(t) \subset  G$ is  a geodesic of  the metric $g^F$ on  $G $. We  may  assume  that
 $g(0) =e$  and $\dot{g}(0) =a \in \gg$.  If $\dot{\g}(0) =w,$  then   the  curve
 $q(t)$   is  a \SR  geodesic if  and only if  its velocity vector  $\dot{q}(0) =
 w + a^*_q  $ is  horizontal,  that is  $\vpi(w + a_q^*) = \vpi(w) +a =0$.
   This proves ii),  which implies  iii). Now iv)
  follows  from  the  remark  that  the  transformation $R_{g(t)}$  deforms  the geodesic   $\g(t)$ in vertical  directions.  Hence  $q(t)$  and  $\g(t)$ have  the  same  projection  to  $M$  (  which   are geodesics   if  and only if   $\g(t)$  is  a  horizontal geodesic.) \qed\\
     Let    $(\pi : Q \to   M, g^Q)$  be  a Riemannian  submersion    and  $D\subset TQ$ a transversal  to   fibers  distribution. Necessary and  sufficient  conditions  when  the projection  to $M$   of $g^Q$-geodesics  coincides  with   projection   of  geodesics of the \SR manifold
      $(Q,D, g|_D)$  are  given   in \cite{M-G}.

  A \SR  geodesic $q(t) =  \tau (\xi(t))$ through a point  $q = q(0)$ is  determined  by   the   initial covector  $\xi(0) \in T_q^*Q$ which may  be  decomposed   as
   $\xi(0)= \xi(0)_{D} + \l $, where   $ \l \in D_q^0 =  \mathrm{Ann}(D)_q$   is  the {\bf codistribution covector}  and   $\xi(0)_D \in  D^*_q$
   is  determined by   the  velocity  vector $\dot{q}(0) \in D_q$.
   The   \SR geodesics,  which  are   horizontal geodesics of  $g^Q$ are  characterized  as geodesics  with  trivial codistribution  covector.
   Comparing   theorem \ref{H-geodesics}   and  theorem \ref{S-geodesics},  we  get
   \bt   Let  $g^D$ be  a  Chaplygin sub-Riemannian   metric  in   a principal bundle
      $(\pi : Q \to  M, \vpi)$ and  $g^Q$   the  standard  extension   of  the sub-Riemannian metric $g^D$.  Then   sub-Riemannian   S-geodesics
       coincide  with H-geodesics  with  trivial codistribution covector.
   \et
\subsubsection{ \label{Bi-invariant case} Bi-invariant  extension of Chaplygin  metric  and Yang-Mills    dynamics}
Assume  now  that $g^Q$ is  a bi-invariant  extension of  the Chaplygin \SR metric
$g^D$,   defined  by  an $\Ad_G$-invariant  metric   $g^{\gg}$ of  the  Lie  algebra $\gg$. Such metric  exists  only  when  $\gg$ is   the  Lie  algebra of  a  compact Lie  group.   Then  the  associated  left  invariant   metric  on  the Lie  group $G$  is  also  right-invariant  and   the   metric  $g^F(a^*, b^*) = g^{\gg}(a,b)  $ on a fiber   $ \pi^{-1}(x)$  and  the   extended Riemannian metric
 $g^Q = g^F + g^D$   are   also  $R_G$-invariant.\\

  In  this  case, the     geodesic  Hamiltonian   system  with Hamiltonian  $h_Q$ has a  nice  physical interpretation   as   dynamical  system,  which  describes    the  evolution  of  a charged   particle in the  base  manifold $M$ in  the  presence of  the Yang-Mills  field,  defined  by   the  principal connection $\vpi :  TQ \to \gg$ ,  see \cite{W},\cite{Mont}. \\
   Recall   that    with  respect  to a   trivialisation
$Q =  M \times G$ ,  the  connection  form  may be  written  as   
 $$    \vpi = \mu + A  = (e^{\a}_L + A^{\a}_i dx^i)\ot e_{\a}  $$
where $(e_{\a})$  is  an orthonormal   basis  of  $\gg$ ,   $(e_{\a}^L)$ (  resp. , $(e^{\a}_L)$)  is   the
corresponding left invariant   field of   frames  (resp.,  coframes)  on  $G$
and  $A = A^{\a}_i(x,g) dx^i \ot e_{\a}$  the  Yang-Mills potential,  given  by  (\ref{connectionform}). The  horizontal  ($R_G$-invariant) lifts  of
the   coordinate  vector  fields  $\p_i := \p_{x^i}$  has  the  form
$\p_i^D := \p_i - A^{\a}_i e_{\a}^R$. Together  with  the fundamental  fields $e_{\a}^*= e^L_{\a}$,
 they form  a frame in $Q = M \times G$.
The  \SR metric is characterized  by  the conditions
$g^D(\p_i^D, \p_j^D) =  g^M(\p_i,\p_j) =  g_{ij}.$
The   vertical metric  $g^F$  is defined    by $g^F(e_{\a}^L, e_{\beta}^L) = g^{\gg}(e_{\a}, e_{\beta}) = \delta_{{\a}{\beta}}$.
  The  metric $g^Q = g^F + g^D$ is  $R_G$-invariant  and   the    fundamental fields $e_{\a}^* =e_{\a}^L$
are  Killing vector  fields.
The  associated   contravariant  metric
 $g^{-1}_Q =  g^{-1}_F + g^{-1}_D$   is  defined  by
 $$g^{-1}_F = \sum e_{\a}^L \ot e_{\a}^L, \,\,\,\,\,\,\,
 g^{-1}_D  = g^{ij}(x) \p_i^D \ot \p_j^D= g^{ij}(\p_i - A^{\a}_i e_{\a}^L)(\p_j- A^{\a}_j e_{\a}^L ) .$$
   Denote  by
 $(x^i, p_i)$  the  local  coordinates in  $T^*M$  with $ T^*M \ni p = p_i dx^i$  and  by  $(g^{\a}, \l_{\a})$  the   local  coordinates  in  $T^*G$,  where $g^{\a}$ are  local  coordinates in  $G$  and $T_g^*G \ni \l = \l_{\a} e^{\a}_L$.  Note  that the linear  forms  $\l_{\a}  \in  T^*_gG$
  are  identified  with  $e_{\a}^L|_g$.  The left invariant  vector  fields
 $\p_i, \p_{p_i}, \p_{\l_{\a}}, e_L^{\a}   $
 form  a  frame  on   $T^*Q = T^*M \times  T^*G $.
 The quadratic  in  momenta Hamiltonians
  $h_M, h_F, h_D , h_Q = h_F + h_D$   can be  written  as follows
$$
\begin{array}{cc}
  h_M =& \frac12 g^{ij}(x)p_ip_j  \\
 h_F =&\frac12 \sum e_{\a}^L e_{\a}^L\\
h_D  =& \frac12 g^{ij}(x)(p_i - A^{\a}_i(x) \l_{\a})(p_j - A^{\beta}_j(x) \l_{\beta}).
\end{array}
 $$
 Using   formula  for  the Poisson  structure  on  $T^*G$, one  can
 easily  calculate  the  Hamiltonian vector  fields and  the  geodesic  equation.  We  consider   another  approach, based on  the O'Neill  formulas.\\

  \bl  The  angle  between a  geodesic $\g(t)$ of  $g^Q$  and  a fundamental field   $a^*, \,  a \in \gg$  is  constant.
  In particular,   the orthogonal projection
  $\mathrm{pr}_{T^vQ} \dot{\g}(t)$  of  the  velocity vector  field  $\dot\g$  to  vertical     subbundle  is  the  restriction  to $\g(t)$  of some fundamental vector  field  $a^*$  and   the velocity  vector  field can be  written  as
   $$   \dot{\g}(t)   = a^*(\g(t) ) + \dot{x}^D(\g(t))$$
  where      $\dot{x}^D(\g(t))$ is    the  horizontal lift  of   the velocity vector   filed  $\dot{x}(t)$ of  the   projection   $x(t)$  of  $\g(t)$ to  $M$.
   \el
  \br Physically,  the   angles  $\varphi_{\a}$ between a geodesic $\dot \g$  and  the  basic   fundamental  fields
$e_{\a}^*$  characterise the  charges of  a particle  with  respect   to   components of  the  Yang-Mills  field  and   the  conditions $\varphi_{\a} =const$ are  called  the  conservation  of  charges. In  particular,  the evolution  of neutral particles  is   described  by  horizontal  geodesics.
\er

\pf
Let  $\g(t)$  be  a  geodesic     and  $ x(t) =
\mathrm{pr}_{M} \g(t)$ its  projection   to  $M$.  Then
  $\dot{x}(t) = \mathrm{pr}_{TM} (\dot {\g}(t))$  and
    the  horizontal    part of   the  velocity vector  field is $\dot{\g}(t)^h = \dot{x}(t)^D$.
   Hence,  we  can  write
    $$\dot{\g}(t) = \dot{x}(t)^D+ u^a(t) e_a^*(\g(t)). $$
  Then
  $$
  \begin{array}{ll}
  \frac{d}{dt}g^Q(e^*_{\beta},\dot{\g}(t) )&= \dot{u}^{\beta}(t)\\
   &= \n_{\dot \gamma}g^Q(e^*_{\beta},\dot{\g}(t) )\\
 & =g^Q(\n_{\dot \gamma}e_{\beta}^*,\dot{\g}(t) )+ g^Q(e_{\beta}^*,\n_{\dot \gamma}\dot{\g}(t))\\
  &=0,
  \end{array}
 $$
 since  the  covariant  derivative  $\n_{\cdot}e_{\beta}^*$  of  a Killing vector  field $e_{\beta}^* = e_{\beta}^L$ is a  skew-symmetric  operator.
  \qed\\

 The  following  theorem   describes  the  relation  between  geodesics of   the Riemannian  metric $g^M$   and     geodesics  of its    bi-invariant  extension $g^Q$.
\bt  A  curve  $\g(t) \subset  Q$  with  projection  $x(t) = \mathrm{pr}_M \g(t)$  and velocity vector  field $\dot{\g}(t) = a^*(\g(t)) + \dot{x}^D (\g(t))$  is a geodesic of  $g^Q$ if  and only  if it satisfies  the equation
  \be \label{Lorentzeq}
  \n^Q_{\dot{\g}(t)} \dot{\g}(t) = (\n^M_{\dot{x}}\dot{x} )^D  + 2\ca^*_{\dot{x}^D}a^* =0.
  \ee
\et

\pf Using O'Neill  formulas,   we  calculate  the  covariant      derivative  $\n^Q_{\dot \g}\dot\g$   of  the   velocity  field
  $\dot \g = a^*(\g(t)) + \dot{x}^D(\dot \g(t))$   as follows
$$ \n^Q_{\dot \g}\dot\g  = (\n^Q_{a^*}a^*)(\g(t)) +\n^Q_{\dot{x}^D}\dot{x}^D + 2 \n^Q_{(\dot{x})^D}a^*  = (\n^M_{\dot x}{\dot x})^D(\g(t)) +
 2  \ca^*_{\dot x^D}a^*(\g(t)). $$
   We use  the fact   that
     geodesics  of  the  bi-invariant  metric  on   a Lie group $G$ are orbits  of
     1-parameter  subgroups,  which implies    $\n^Q_{a^*}a^* = \n^F_{a^*}{a^*} =0$.
 \qed

 Recall  that    $2\ca^*_{X^D} a^* = - F^*_X \l , \,  \l =  g^Q \circ a^* \in  D^0$
 where   $F^*_X : \gg^* \to   \G D  $ is   the  linear  map,  dual  to  the
 map $F_{X_q} : D_q \to   \gg$,  associated  with   the   curvature  2-form $F$.

 The equation (\ref{Lorentzeq}) is   equivalent  to  the   equation
\be \label{Lorentz1}
  \n^M_{\dot x} \dot x  = g_M^{-1} \l_bF^b_i(\dot x)
 \ee

 where the   right hand  side  is  the vector  field   metrically  dual   to
  the   1-form   $  \l_bF^b_i(\dot x) $   and   $\l \in \gg^*$ is  a  constant  covector   (a  charge). \qed\\

  The   equation (\ref{Lorentz1})  describes  the  motion   of   a charged particle  in   the   Yang-Mills  field    $\vpi$
    with    the strength  tensor  $F$.
   In  the  case  when $\pi : Q \to M$ is  a  circle bundle,   the  connection
    $\vpi$   defines  the  Maxwell   field   (if   $g^M$  has Lorentz  signature)    and  the  equation reduces  to  the  Lorentz  equation
     for  a  charge  particle  in   the electromagnetic  field,  defined  by  the curvature  2-form  $F$.

\section{Homogeneous  sub-Riemannian   manifolds }
We consider  some  class of  homogeneous  \SR manifolds, for  which  S-geodesics  coincide  with  H-geodesics and describe  sub-Riemannian  symmetric  spaces.\\

\subsection{  Chaplygin    system on  homogeneous  spaces}
\subsubsection{Chaplygin  system of  a  Lie  group}
    Let  $\pi :  G \to   M= G/H$    be    the principal bundle  associated  to  a homogeneous
     Riemannian manifold  $(M = G/H, g)$. A   reductive  decomposition  $\gg =\gh  + \gm,\,  \gm = T_oM, \, o =eH$,  defines
     a principal  connection with    connection  form
     $\vpi  = \mathrm{pr}_{\gh} \circ \mu^L$,   which is   the  projection to  $\gh$ of  the left invariant Maurer-Cartan  form  $\mu^L$.
     Denote  by  $(D = \ker \vpi, g^D)$ the associated Chaplygin \SR  metric. Since the  stability  subalgebra $\gh$ is  compact, it admits a bi-invariant  Euclidean metric $g^{\gh}$. We denote  by $g^{G}$ the associated  bi-invariant  extension of the \SR metric. It is  a left  $G$-invariant and  right $H$-invariant  metric   on  $G$. \\
     The     distribution  $D : =\ker \vpi $  with    $D_0 = \gm$
      is bracket generating if and only if  $\gm$ generates   the Lie  algebra  $\gg$. The  Jacobi identity  shows  that
       $\gg': = [\gm,\gm] + \gm$ is  a  subalgebra of $\gg$. It  generates  a  subgroup $G' \subset G$ which  acts transitively in $M$. Hence,
         changing   $G$ to $G'$,  we   may always  assume  that  $D$ is bracket generating.

    \bp  Let   $(M= G/H, g^M)$  be  a  homogeneous Riemannian manifold  with   reductive  decomposition
     $\gg = \gh + \gm $  such  that $\gm$ generates $\gg$. Then the principal  connection
     $\vpi:= \mathrm{pr}_{\gh} \circ \mu^L $  defines  a  Chaplygin left invariant  sub-Riemannian  structure
      $(G, D, g^D)$ on the Lie  group $G$  with  the  connection $\vpi$. A bi-invariant  metric $g^{\gh}$  on  $\gh$ defines
       a bi-invariant extension  of  the sub-Riemannuian  metric  $g^D$ to a left invariant    Riemannian metric   on  $G$.
     \ep

 \subsubsection{ Chaplygin  systems on  homogeneous  manifolds}
     Now  we  consider  a  generalisation   of  the  above  construction.

Assume  that  the  stabilizer $H$ of  a Riemannian  homogeneous manifold  $M = G/H$ is  an  almost   direct product
$H = K \cdot  L$ of   two  compact normal  subgroups  and   $\gk, \gl$   are  associated  Lie  subalgebras.
Then   $\pi : Q = G/K \to  M =  G/K \cdot L$ is an $L$-principal  bundle   with   the right   action  of  $L$  and
  $M = G/H$  has    the   reductive  (i.e.  $\Ad_H$-invariant )  decomposition   of    the  form
 $$\gg = \gh + \gm  = (\gk  + \gl) + \gm.$$
  The projection    $\vpi^G = \mathrm{pr}_{\gl} \mu^L:   TG \to \gl$
 of the left invariant Maurer-Cartan  form  $\mu$  to $\gl$ is     a  left invariant   $\gl$-valued 1-form.
    The  form  $\vpi^G$ is  right $K$-invariant 
     and  right   $L$-equivariant,  that is
      $R_\ell^* \vpi = \Ad^{-1}_\ell \circ \vpi , \,  \ell \in L $.   Hence  it projects  to  a  $G$-invariant principal  connection  form
       $\vpi : TQ  \to  \gl$.   The principal  bundle  $\pi : Q= G/K \to  M = G/K \cdot L$ with  the  connection  form
        $\vpi : TQ \to \gl$  defines a Chaplygin \SR metric  $(D = \ker{\vpi}, g^D)$.  As  above,  it  admits  a bi-invariant  extension. We  get

       \bp A homogeneous Riemannian manifold   $(M = G/H , g^M)$ with non  simple  stabilizer  $H =  K \cdot L$
  defines  an invariant sub-Riemannian  Chaplygin metric  $(D, g^D)$  on  the  total  space  of  the principal $L$-bundle  $\pi : Q= G/K  \to  M =  G/K \cdot L$  with the connection  form  $\vpi : TQ \to  \gl$,  which is  the projection to  $Q$ of  the  form  $\vpi^G = \mathrm{pr}_{\gl}\circ \mu^L$ on $G$. The  sub-Riemannian  metric  admits    a bi-invariant  extension  to
  an invariant metric   $g^Q$ on  $Q$.
  \ep
\subsubsection{Homogeneous  contact     sub-Riemannian   manifolds}
   The  above  construction  may be  applied   to  homogeneous  Sasaki manifolds.        
   We consider the
   case   of   regular  compact  homogeneous  Sasaki manifolds,  described   as  follows.
     Let  $M = G/H$     be  a  flag manifold  (i.e.  an  adjoint orbit  of  a  compact  semisimple Lie  group $G$)
       and $g = \o(\cdot, J \cdot)$  an  invariant Hodge-K\"ahler  metric on  $M$,  where  $ J$
      is  an invariant complex  structure and $\o $  an integer invariant  symplectic  form    (the K\"ahler form).
   Then   there  exists   a  homogeneous principal circle  bundle  $ \pi : Q =  G/K \to M =  G/H =  G/K \cdot S^1$
   with  a  principal  connection   $\vpi : TQ \to \bR = \mathrm{Lie{S^1}}$, whose  curvature  form is $\o$.
   The  K\"ahler metric  is naturally extended  to  an invariant Sasaki metric  $g^Q$,  such  that  the  fundamental  field
    $Z$  of  the $S^1$-bundle $\pi$   is a Killing  field. \\
   This   Sasaki metric  $g^Q$   is   the bi-invariant  extension  of  the Chaplygin  \SR metric  $(D,g^D)$ associated  to  the principal
    bundle  $\pi : Q \to M$ with  the   connection $\vpi$.\\

From physical point of  view,  the  principal $S^1$-bundle  $\pi :  Q \to M$ with Sasaki metric   corresponds to Kaluza-Klein  description  of electromagnetic  field.
    The  projection  to  $M$ of   geodesics of Sasaki metric  are  solutions  of  the Lorentz  equation  which  describes   the  evolution  of electric  charges in the electromagnetic  field $\o$.

\subsection{Symmetric  sub-Riemannian   manifolds}

 Strichartz  \cite{Str}     defined  the  notion  of {\bf \SR symmetric  space}
as  a homogeneous \SR  manifold   $(Q = G/H,D, g^D)$  such  that  the  stabilizer $H$  contains  an  involutive  element $\s$ (called the sub-Riemannian   symmetry)  which  acts on  the   subspace   $D_o$ at  the point   $o =  eH  \in Q$    as  $- \mathrm{id}$.\\
 He   classified  3-dimensional  \SR  symmetric  spaces and   stated  the problem of  extension of  this classification    to higher  dimensions.\\
P. Bieliavsky, E. Falbel  and  C. Gorodski  \cite{B-F-G}    classified   symmetric  \SR manifolds of  contact  type.
W. Respondek  and A.J. Maciejewski \cite{R-M} describe  all integrable \SR metrics on 3-dimensional Lie groups  with integrable  H-geodesic  flow. They  are  exhausted by \SR   symmetric  spaces. \\
Below   we  recall  basic properties of  affine  symmetric  spaces and   give a construction  of   \SR  symmetric  spaces   in terms of   affine  symmetric  spaces: Any bracket generating  \SR  symmetric  space is  the  total  space $M = G/K$ of  a homogeneous  bundle   $\pi : M = G/K \to  S = G/H$  over    an   affine  symmetric  space  $S = G/H$, determined  by  a  compact  subgroup  $K$ of  the  stability group  $H$.

\subsubsection{Affine  symmetric  spaces}
Let $(M, \n)$  be a  (connected) manifold  with a linear  connection $\n$.  A non-trivial   involutive  automorphism $\s=\s_x$ of  $(M,\n))$
is called  a {\bf cental  symmetry  with center $x \in M$} if $\s$ preserves $x$ and  acts  as $- \id$ in  the tangent  space  $T_xM$. The manifold $(M,\n)$ is called  an {\bf (affine) symmetric  space } if  any point  is  the center of  some central  symmetry  $ \s_x$.  A product  $\s_x \s_y$ of  two central  symmetries   with  sufficiently  closed  to each other  centers $x,y$ is  a  shift  along  the  geodesics, connecting  these points. This  implies  that the   group $G$, generated  by  all central symmetries  is  a transitive  Lie group, called  the {\bf  transvection  group}. The manifold  $M$
 is identified  with  the quotient  space  $M = G/H$,  where  $H$ is   the stabilizer of  a point  $o \in M$. Then the  central  symmetry $\s=\s_o$  defines  an involutive automorphism
$s = \Ad_{\s_0} : g \mapsto s(g):=\s_0\circ g \circ \s_0$  of the Lie  group   $G$, which acts  trivially on the  connected component $H^0$   of $H$.
 We  denote     by $s$  also the induced involutive  automorphism of  the  Lie  algebra   $\gg = Lie (G) $. Its  eigenspace  decomposition
$$\gg =  \gg_+ + \gg_-  , \,\,  s|_{\gg_{\pm}}  = \pm \id,   $$
  where  $\gg_+ = \gh = \mathrm{Lie}(H)$,  is called  the {\bf  symmetric  decomposition}.
It is characterized     by  the  conditions
 $$   [\gg_+,\gg_-] \subset \gg_-,\,\,\, [\gg_-,\gg_-] \subset \gg_+.$$
  Moreover, if $G$ is  the transvection  group,  then
  \be \label{transvectionproperty}
  [\gg_-, \gg_-] = \gg_+.
  \ee

  The   geodesics through   the point  $ o = eH$   are orbits   $e^{tX}o$ of  1-parametric  subgroups  $e^{tX} \in  G$ generated  by  elements  $X \in \gg_-$.

The     following well known  result establishes  a bijection between   symmetric  decomposition $\gg  = \gg_+ + \gg_- $ of  a  Lie  algebra $\gg$   with (\ref{transvectionproperty}) and simply  connected  affine  symmetric  spaces  $S=G/H$,  where $G$  is  the
   simply  connected transvection  group  with  $\mathrm{Lie}(G)= \gg$.
\bt  Let      $\gg  = \gg_+ + \gg_- $ be  a  symmetric  decomposition with (\ref{transvectionproperty})
 associated to an involutive   automorphism  $s$.    Denote  by   $G$ the  simply  connected  Lie  group  with     $\mathrm{Lie}( G)=\gg$
   and by
 $H^0$ the  connected  subgroup of $G$, generated  by $\gh = \gg_+$. Then   $S = G/H^0$ is  a simply  connected  affine  symmetric  space.
  The  invariant torsion  free  linear  connection  $\n$ in  $S$ is  defined  by the  condition
$$  \n_{X}Y^*|_o = -\frac12 [X,Y]_o, \,\,  X,Y \in  \gg_-  =  T_oS  $$
where   $Y^*$ denote the velocity vector  field   of a 1-parameter subgroup $e^{tY},\,  Y \in \gg_-$.
The  central  symmetry   with the center $o = eH^0$ is defined  by \\
$$\s :  x = gH^0 \mapsto  \s x := s(g)H^0$$
where $s$ is the involutive  automorphism of the Lie group $G$  
 generated  by the  automorphism  $s$  of  $\gg$.
 Moreover,  {\bf any}  affine  symmetric  space, associated  with  the  above symmetric  decomposition,  has  the  form  $G/H$
 where  $H$  is  a closed subgroup such  that   $H^0 \subset H \subset G^{\s}$. Here $G^{\s}$ is  the  fixed point  set of
 $\s$.
\et

\subsubsection{Sub-Riemannian   symmetric spaces    associated  with an  affine  symmetric space}

  Let  $(S= G/H,  \n, \s)$   be  a simply  connected   affine  symmetric   space with   the  transvection  group $G$.Without loss of generality, we may  assume  that       the central  symmetry
  $\s = \s_o$  belongs  to  the center  $Z(H)$ of  the  stability  subgroup.Then the  associated involutive  automorphism
  $s = \Ad_{\s}$  of $G$  acts  trivially on  $H$ and  defines
        a symmetric  decomposition  $\gg = \gg_+ + \gg_-$  where $\gg_+ = \gh = Lie (H)$.
  Let $K \subset H$ be  a  compact  subgroup of  $H$ which  contains $\s$.\\
  The  homogeneous manifold  $Q = G/K$ has   a reductive  decomposition
    $$  \gg= \gk + \gm =  \gk + (\gp + \gg_-)  $$
      where $\gk = Lie (K)$  and $\gg_+ = \gk + \gp$ is  a reductive  (i.e.  $\Ad_K$-invariant )
  decomposition  of $\gg_+$. We identify $\gm $ with  the tangent  space  $T_o Q$ at  the point $o =eK$.  Then  the isotropy  representation  of $ K $ in
  $T_oQ$ is identified  with  $ \Ad_K|_{\gm}$.
   The  $\Ad_K$-invariant  subspace $\gg_- $ of  the tangent   space $\gm = T_oQ$  is naturally  extended  to an invariant  distribution $D \subset TQ$. More precisely, for   $x = aK \in G/K$,  the subspace
   $D_x =( L_{a})_* \gg_-$, where   $L_a : bK \to abK$ is  the  action  of $G$ in  $Q=G/K$.\\
     The  distribution $D$ is invariant  with  respect to the action of involution  $\s \in  K $ and  the  isotropy  action  $\Ad_{\s}|_{\gm}$
     of $\s$ acts   on  $\gg_- \subset \gm$  as $ - \id$.  Since $G$ is   the transvection  group, $[\gg_-, \gg_-] = \gg_+$  and $D$ is bracket generating distribution.
     Since  the  group $K$ is  compact,  there exist  an $\Ad_K$-invariant Euclidean metric $g$ in $\gg_-.$ It is naturally  extended  to  an  invariant
     \SR metric $g^D$ in  $D$, defined  by
     $$   g_x^D(X,Y) := g( (L_a^{-1})_* X, (L_a^{-1})_*Y  )   ,\,\,\,  a \in G, \, x = aK, \, X,Y \in  T_xQ = (L_a)_* \gg_-.                $$
      Hence, the invariant \SR manifold $(Q= G/K, D, \s)$  is a \SR  symmetric  space.
     This  proves  the first claim of  the  following  theorem.

\bt i)  Let   $(S = G/H, \n, \s)$   be    a simply  connected   affine   symmetric  space with  the transvection  group $G$   
     and   $K \subset H$  a compact  Lie subgroup, which  contains $\s $ as  a central element.   Let  $g$  be an  $\Ad_K$-invariant Euclidean  metric    in  $\gg_-$. Then the Euclidean  space
     $(\gg_-, g)$ is  extended  to an invariant sub-Riemannian  structure $(D, g^D)$ in  $Q =  G/K$    such  that $(Q=G/K,D, g^D, \s)$ is  a bracket generating sub-Riemannian  symmetric  space.

ii) Conversely, up to a  covering any bracket generating   sub-Riemannian symmetric  space  can be obtained  by  this  construction.
\et
\pf   ii)  Let   $(Q= G/K, D, g^D, \s)$ be  a bracket generating \SR symmetric  space,  where $\s \in K$ is    the \SR  symmetry with center $o = eK$   and    $s = \Ad_{\s}$ the  associated involutive  automorphism  of  $G$ and  $\gg$.
 We   may  chose  a reductive  decomposition
 $$\gg= \gk + \gm  =   (\gk_+ + \gk_-) + (\gm_+ + \gm_-)$$ of  $G/K$ which is  consistent  with  the  symmetric  decomposition  $\gg = \gg_+ + \gg_-$, defined  by  $s$, such  that  $\gg_+ = \gk_+ + \gm_+,\,  \gg_- = \gk_- + \gm_-$. By definition,the subspace $D_o = D|_o \subset T_oQ = \gm$  belongs  to  $\gm_-$. Since  the  distribution $D$ is bracket generating, we may  assume  that  the  subalgebra  $\bar{\gg}$ generated by $\D_o$ coincides  with $\gg$.  But $\bar{\gg} = [D_o, D_o] + D_o\subset \gg_+ + D_0$. This implies  that  $\gg= \gk_+ + \gm_+ + \gg_-$, and   $\gk= \gk_+ , D_o = \gg_-, \gm = \gm_+ + \gg_- = \gg_-$. Denote  by   $G_+$  the connected subgroup  of $G$, generated by $\gg_+$, Since it  commutes  with $\s$, it is   the  connected  component   of  the  group $ H= G_+ \cup \sigma G_+$. The manifold $S =  G/H$ is  an  affine  symmetric  space    with the symmetry $\s \in H$, belonging  to the center. Consider  the  subgroup   $K' = K \cap H$ with  the Lie  algebra $\gk$. It also contains $\s$  as  a central  element. The claim i) shows that the  space  $Q'= G/K'$  has  a structure  of  \SR  symmetric  space  which is  locally isomorphic  to  the  initial  \SR  symmetric  space   $Q$. \qed

\subsubsection{Compact sub-Riemannian  symmetric  space  associated  to a  graded   complex  semisimple Lie   algebra}
We show   that   any flag manifold of depth $>1$ admits    the   structure  of  bracket generating symmetric   \SR  manifold.

Let   $\gg = \sum_{i=-d}^d\gg_i$
 be  a  {\bf fundamentally graded } complex  semisimple  Lie  algebra  of  depth  $d \geq 2$ ( s.t.     $\gg_{-1}$  generates  $\gg_-$) and
 $\gp := \sum_{i\geq 0} \gg_i$ the  associated parabolic  subalgebra.  The  associated ( complex  compact simply  connected) homogeneous manifold $F= G/P$,  where  $G\supset P$ are the Lie  groups  associated to Lie  algebras $\gg \supset \gp$, is called  a {\bf  flag manifold}. \\
 Denote  by  $\tau$ the   anti-linear involution of $\gg$, which  defines  the  compact  real form $\gg^{\tau}$  s.t.
   $\gg^{\tau}  = \gg^{\tau}_0  + \sum_{i>0} \gm_i ,\,\, \gm_i:= (\gg_{-i} +\gg_i )^{\tau}.$\\
    The Lie algebra  $\gg^{\tau}$  has  the  symmetric  decomposition
  $$\gg^{\tau}   =  \gg^{\tau}_{ev} + \gg^{\tau}_{odd}=(\gg^{\tau}_0 +
   \sum_{i\equiv 0(\mathrm{mod}2)} \gm_i ) + \sum_{i \equiv 1 (\mathrm{mod} 2)}\gm_i .  $$
    We denote  by $s$  the  associated  involution  of  the Lie  algebra $\gg^{\tau}$ and the  corresponding   simply  connected compact  Lie  group   $G^{\tau}$.
 Denote  by
 $H \subset G$ the connected  compact  subgroup generated  by  $\gh = \gg_0^{\tau}$.
 The  group  $G^{\tau}$ acts transitively on  the  flag  manifold $F$ with  stability subgroup  $H$   and   has the reductive  decomposition
  $$\gg^{\tau} = \gh + \gm = \gg_0^{\tau} + \sum_{i>0}\gm_i.$$
 The  involutive  automorphism $s$ acts by  $\s|{\gh} = \id,\,\, \s|_{\gm_i} = (-1)^i \id.$\\
  Denote  by  $D \subset TF$ the   ( bracket generating) invariant  distribution generated by
   $\gm_{1}$ and  by $g^D$ the invariant  \SR metric in $D$  defined  by an $\Ad_H$-invariant metric in $\gm_{1}$.  Then   $(D, g^D)$ is  an invariant \SR metric of $F = G/H$. Moreover,
   $(F= G^{\tau}/H, D, g^D)$ is  a \SR   symmetric  space, where  the  symmetry, defined  by  the involutive  automorphism  $s$.

    This implies
    \bt   Let  $\gg = \sum_{i=-k}^k \gg_i$ be a  fundamental  depth  $k>1$ gradation  of a  complex semisimple Lie  algebra
     and let $F = G/P$ be  the associated   flag manifold. Denote by
    $$   \gg^{\tau} = \gh + \gm = \gh + \sum_{i=1}^k \gm_i, \,\gh = \gg_0^{\tau,\,} \gm_i = (\gg_{-i} + \gg_i)^{\tau}   $$
     the  associated  decomposition  of  the   compact real form $\gg^{\tau}$ and by
     $  g^{\gm_{1}}$  an $\ad_{\gh}$-invariant Euclidean metric in $\gm_{1}$. \\
     Then  the pair $(\gm_1, g^{\gm_1})$ defines  an invariant bracket generating  sub-Riemannian  metric $(D, g^D)$ on  the  flag manifold    $F=G^{\tau}/H $
     considered  as  a homogeneous manifold  of  the  compact real  form  $G^{\tau}$  of  $G$.
 Moreover, the  sub-Riemannian manifold  $(F=G^{\tau}/H, D, g^D)$ is  a sub-Riemannian  symmetric  space   with  the  symmetry  defined  by the involutive  automorphism $s$ of $\gg$,  associated with  the  symmetric  decomposition $\gg^{\tau} = \gg^{\tau}_{ev} + \gg^{\tau}_{odd}$.
 \et

{ \bf Example}
    Let
    $$\gg  = \gg_{-2} + \gg_{-1} + \gg_9 + \gg_0 + \gg_1 + \gg_2,\, \dim  \gg_{\pm 2} =1$$
  be    the  contact gradation of  a  complex  simple  Lie  algebra $\gg$ , i.e.  the  eigenspace  decomposition  of $\ad_{H_{\mu}}$  where  $H_{\mu}$ is  the  coroot   associated  to    the maximal  root   $\mu$   of  $\gg$.  Then   the   symmetric  space  $G^{\tau}/ G^{\tau}_{ev}$  is  the   quaternionic  K\"ahler  symmetric  space  ( the Wolf  space )  and   the  flag manifold $F = G^{\tau}/H$, where $  Lie(H) = \gh_0^{\tau}$,  is  the  associated  twistor  space.   The distribution  $D$ is  the holomorphic  contact distribution  and   $g^D$ is the unique (up  to scaling )  invariant \SR metric  on  $D$ ( for  $\gg \neq \gsl_n(\bC)$ ).  It is  the  restriction  of   the invariant K\"ahler-Einstein metric  on $F$.

\end{document}